\documentclass[11pt,a4paper,reqno]{amsart}
\usepackage{amsmath,amsfonts,amsthm,amssymb,graphicx,amscd,color,colortbl}
\usepackage{float,epsf,subfig}
\usepackage[mathscr]{euscript}
\usepackage[english]{babel}

\usepackage[latin1]{inputenc}
\usepackage{verbatim}
\newenvironment{my_itemize}{
\begin{itemize}
  \setlength{\itemsep}{2mm}
  \setlength{\parskip}{0pt}
  \setlength{\parsep}{0pt}}{
  \end{itemize}
}

\theoremstyle{definition}

\theoremstyle{remark}

\numberwithin{figure}{section}
\numberwithin{equation}{section} 

\newcommand{\ci}{\cite}
\newcommand{\bib}{\bibitem}
\newcommand{\uu}{\mathbf{u}}
\newcommand{\Ka}{\mathcal{K}}

\providecommand{\abs}[1]{\lvert#1\rvert}
\providecommand{\norm}[1]{\lVert#1\rVert}

\makeatletter
\newcommand{\rmnum}[1]{\romannumeral #1}
\newcommand{\Rmnum}[1]{\expandafter\@slowromancap\romannumeral#1@}
\makeatother

\begin{document}
\title[$L^1$--Stability of Vortex Sheets and Entropy Waves]{$L^1$--Stability of Vortex Sheets
and Entropy Waves\\ in Steady Supersonic Euler Flows over Lipschitz Walls}
\author{Gui-Qiang G. Chen $\qquad$ Vaibhav Kukreja}

\address{Gui-Qiang G. Chen, Mathematical Institute, University of Oxford,
         Oxford, OX1 3LB, UK}
 \email{\tt chengq@maths.ox.ac.uk}

\address{Vaibhav Kukreja, Moshman Research, Portland, OR 97219, USA; Department of Mathematics, Northwestern University,
         Evanston, IL 60208, USA}
\email{vaibhavk@moshmanresearch.com}




\date{\today}
\maketitle

\begin{abstract}
We study the well-posedness of compressible vortex sheets and entropy waves
in two-dimensional steady supersonic Euler flows over Lipschitz walls with $BV$ incoming flows.
Both the Lipschitz wall of $BV$ tangential angle function
and the $BV$ incoming flow perturb a background
strong vortex sheet/entropy wave.
In particular, when the total variation of the incoming flow perturbation
around the background strong vortex sheet/entropy wave is small,
we prove that the two-dimensional steady supersonic Euler flows
containing a strong vortex sheet/entropy wave past the Lipschitz wall are $L^{1}$--stable.
The weak waves are reflected after the nonlinear waves interact
with the strong vortex sheet/entropy wave and the wall boundary.
Using the wave-front tracking method, the existence of solutions in $BV$ over the Lipschitz walls
is first shown, when
the total variation of the incoming flow perturbation around the background strong vortex sheet/entropy
wave is suitably small.
Then we establish the $L^{1}$--stability
of the solutions with respect to the incoming flows.
To achieve this, a Lyapunov functional,
equivalent to the $L^{1}$--distance between two solutions containing
the strong vortex sheets/entropy waves, is carefully constructed to
include the nonlinear waves generated by both the wall boundary and the incoming flow.
This Lyapunov functional is then proved to decrease in the flow direction, leading to
the $L^{1}$--stability of the solutions.
Furthermore, the uniqueness of these solutions extends to a larger class of viscosity solutions.
\end{abstract}
\maketitle

\bigskip
\medskip

\section{Introduction}
We study the well-posedness of
compressible vortex sheets and entropy waves
in two-dimensional ($2$-D) steady supersonic Euler flows over Lipschitz
walls with $BV$ incoming flows.
The inviscid compressible flows are governed
by the $2$-D steady Euler system:
\begin{equation}\label{1.1}
\begin{cases}
(\rho u)_x + (\rho v)_y = 0, \\
(\rho u^2 + p)_x + (\rho uv)_y = 0,  \\
(\rho uv)_x + (\rho v^2 + p)_y = 0,  \\
(\rho u(E + \frac{p}{\rho}))_x + (\rho v(E + \frac{p}{\rho}))_y = 0,
\end{cases}
\end{equation}
with $\uu=(u,v)$, $p$, $\rho$, and $E$ representing the fluid velocity, scalar pressure, density,
and total energy, respectively. Furthermore, the total energy $E$ is explicitly given by
\begin{equation*}
E = \frac{1}{2}|\uu|^2 + e,
\end{equation*}
where the internal energy $e$ can be written as a function of $(p,\rho)$
defined through the thermodynamical relations:  $e=e(p,\rho)$.
The other two thermodynamic variables are  entropy $S$ and temperature $T$.
In the case of an ideal gas, pressure {\it p} and internal energy {\it e} can be expressed as
\begin{equation}\label{1.2}
p = R \rho T, \hspace{7mm} e = c_{\nu}T
\end{equation}
with the adiabatic exponent $\gamma= 1 + \frac{R}{c_{\nu}} > 1$.
In particular, in terms of density $\rho$ and entropy {\it S}, they have the form:
\begin{equation}\label{1.4}
p(\rho, S)= \kappa\rho^{\gamma}e^{S/c_{\nu}}, \qquad
e(\rho,S)= \frac{\kappa}{\gamma - 1}\rho^{\gamma - 1}e^{S/c_{\nu}}=\frac{R}{\gamma-1}T(\rho,S),
\end{equation}
where $R, c_{\nu}$, and $\kappa$ are all positive constants.

When entropy $S$ is constant, the flow becomes {\it isentropic}, which is governed
by the first three equations in \eqref{1.1} with
the pressure-density relation as $p=p(\rho) = \frac{\rho^{\gamma}}{\gamma}$ for $\gamma>1$.
The limiting case $\gamma = 1$ corresponds to the isothermal flow.

Define
\begin{equation*}
c := \sqrt{p_{\rho}(\rho, S)}  \nonumber
\end{equation*}
as the sonic speed. Then, for polytropic gases,
the sonic speed is $c = \sqrt{\gamma p/\rho}$.
The type of flow is classified by the {\it Mach number} $M:= \frac{|\uu|}{c}=\frac{\sqrt{u^2+v^2}}{c}$.
When $M > 1$, system \eqref{1.1} governs a {\it supersonic} flow
({\it i.e.}, $|\uu| > c$), which has all real eigenvalues and is hyperbolic.
For $M < 1$, system \eqref{1.1} governs a {\it subsonic} flow
({\it i.e.}, $|\uu|< c$), which has complex eigenvalues and
is mixed-composite elliptic-hyperbolic.
A {\it sonic} state corresponds $M = 1$.

We are interested in the $L^1$--stability of compressible vortex sheets and entropy waves in steady supersonic flow
over the Lipschitz walls under the $BV$ perturbations of the incoming
flow (see Fig. 1.1).
Multidimensional (M-D) steady supersonic Euler flows are important in many physical applications
({\it cf.} Courant-Friedrichs \ci{Courant-Friedrichs-1948}).
In particular, when the upstream flow is a uniform steady
flow above the plane wall in $x < 0$ all the time,
the downstream flow above a Lipschitz wall in $x > 0$ is governed by a steady Euler flow
after a sufficiently long time.
Moreover, compressible vortex sheets and entropy waves occur ubiquitously in nature and
are fundamental waves.
Furthermore, since steady Euler flows are large-time asymptotic states and may be global attractors
of the corresponding unsteady compressible Euler flows,
it is important to establish the existence of
steady Euler flows and understand their qualitative properties,
which are still wide open.
In particular, the uniqueness and stability of compressible vortex sheet/entropy wave solutions
in a class of physical entropy solutions for steady supersonic flow
has been a longstanding open problem.
On the other hand,
compressible vortex sheets and entropy waves may be formulated as characteristic free boundaries,
and the stability problem can also be formulated as a free boundary problem ({\it cf}. \cite{Chen-Feldman2018,CSV,CW12}),
whose solution is a direct corollary of the well-posedness results in $BV$ established in this paper,
for which the regularity of the free boundary in $\mathbb{R}^2$ is Lipschitz (with bounded total variation of the tangential angle function)
but not in $C^1$ in general.

The stability of contact discontinuities for the Cauchy problem for
$1$-D strictly hyperbolic systems under a $BV$ perturbation has been
studied by Sabl$\acute{\text{e}}$-Tougeron \ci{Sable-1993}
and Corli--Sabl$\acute{\text{e}}$-Tougeron \ci{Corli-Sable-1997}.
In particular, when a weak wave interacts with the boundary of
strip $\{(t,x)\,:\, t \ge 0, -1 < x < 1\}$,
the reflection coefficients for the reflected waves (similar to $K_{11}$ in \eqref{key-1} below) are required to be
less than one, which is the stability condition for the mixed problem in the strip in
the earlier works; see, {\it e.g.}, Sabl$\acute{\text{e}}$-Tougeron \ci{Sable-1993}.
The nonlinear structural stability with the local-in time existence of $2$-D compressible vortex sheet solutions was first established
for the Mach number $M>\sqrt{2}$ in Coulombel-Secchi \cite{CS1,CS2},
while the $2$-D compressible vortex sheets are not stable in general even locally when $M<\sqrt{2}$; also see \cite{CW12} for further results
for compressible vortex sheets.
Moreover, multiple wild solutions for the Cauchy problem of the
compressible Euler equations have been constructed; see \cite{CDK,KKMM} and the references cited therein
for both the isentropic and full Euler cases.
Thus, it is fundamental to understand further the underlying physics of the stability/instability of
compressible vortex sheets and entropy waves
and their interactions with other nonlinear waves,
even for the large Mach number cases.
In particular, it is important to understand whether strong steady compressible vortex sheets and entropy waves
are $L^1$-stable in the class of entropy solutions in $BV$
for steady supersonic flow for any Mach number $M>1$, different from the time-dependent case in \cite{CS1,CS2} as the results of this paper have indicated.
We hope that the analysis and results in this paper will inspire further physics-based modeling and analysis of the interactions between
the two types of strong waves and other nonlinear waves,
and possibly light a path forward to further understanding of the global existence and nonlinear stability of vortex sheet/entropy wave solutions
for the $2$-D compressible Euler equations in gas dynamics.

Working with the full Euler system \eqref{1.1} and a uniform upstream flow containing
one straight strong vortex sheet/entropy wave,
Chen-Zhang-Zhu \ci{Chen-Zhang-Zhu-2007}
first established the global existence of supersonic Euler flows in $BV$ with
a strong vortex sheet/entropy wave under the $BV$ perturbation of the Lipschitz wall
by using the Glimm scheme.
The essential difference between system \eqref{1.1} (as analyzed in \ci{Chen-Zhang-Zhu-2007}
and \S 2--\S 7 below) and strictly hyperbolic systems as
considered in \ci{Corli-Sable-1997, Sable-1993}
is that two of the four characteristic eigenvalues coincide and have two corresponding
linearly independent
eigenvectors that determine precisely the compressible vortex sheet and entropy wave,
so that two independent parameters are required to describe them, respectively.

\begin{figure}[tbh]
\centering
\includegraphics[width =10 cm, height =5.5 cm] {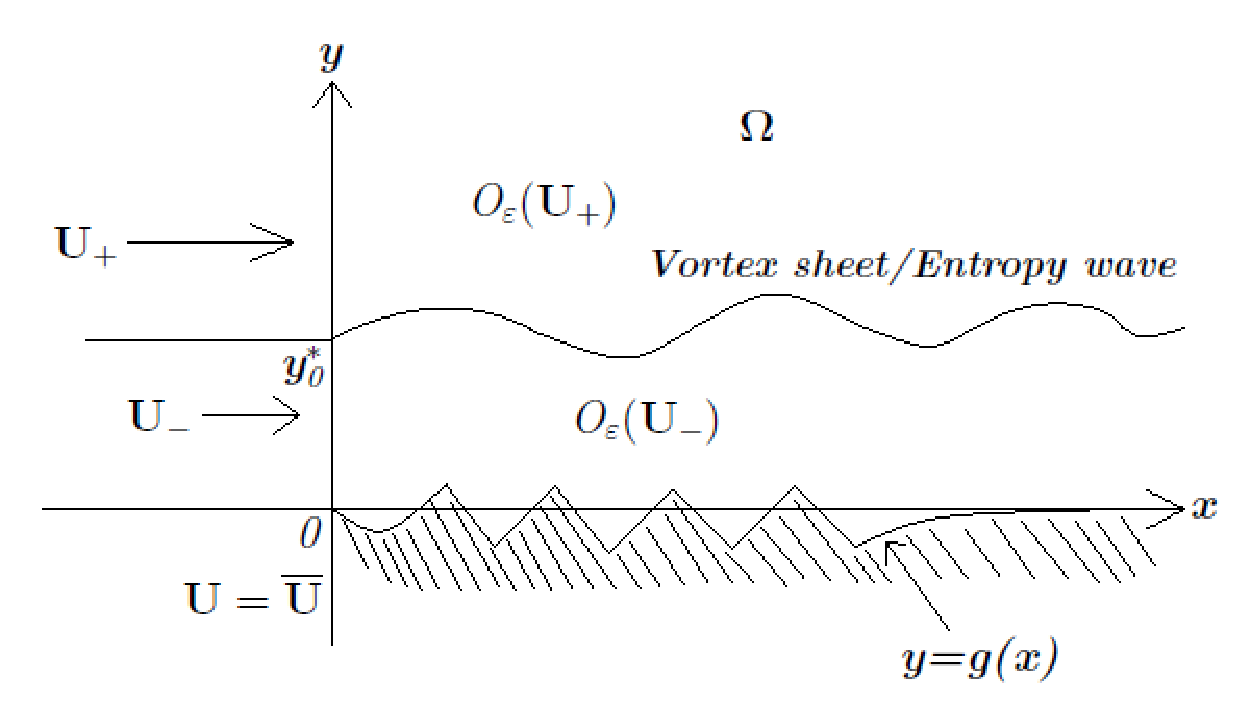}
\caption{Stability of the compressible vortex sheet/entropy wave in supersonic flow} \label{Figure1}
\end{figure}

Consider the following vector functions of the physical variables $U=(\uu, p,\rho)^\top=(u,v,p,\rho)^\top$:
$$
W(U) = (\rho u, \rho u^2 + p, \rho uv, \rho u(h + \frac{u^2 + v^2}{2}))^\top,
$$
$$
H(U) = (\rho v, \rho uv, \rho v^2 + p, \rho v(h + \frac{u^2 + v^2}{2}))^\top,
$$
with $h = \frac{\gamma p}{(\gamma - 1)\rho}$.
Then the steady Euler equations in \eqref{1.1} can be expressed in the
following conservative form:
\begin{equation}\label{2.1}
W(U)_x + H(U)_y = 0.
\end{equation}
In this paper, for completeness,
we first show, via the wave-front tracking method, the existence of solutions
of the problem when a small $BV$ perturbation is added to the uniform incoming flow
containing one straight strong vortex sheet/entropy wave.
Then the $L^1$--stability of entropy solutions containing strong vortex sheets/entropy waves
is established.
As corollaries of these results, the estimates on the uniformly Lipschitz semigroup
of entropy solutions generated by the wave-front tracking approximations are obtained,
and the uniqueness of solutions containing strong vortex sheets/entropy waves
is established in a larger class of solutions, {\it i.e.},
the class of viscosity solutions.
More precisely, we focus mainly on the problem in domain $\Omega$ over the Lipschitz wall
for the supersonic Euler flows $U$ governed by system \eqref{2.1}, given that
the corresponding problem for the isentropic system is simpler to analyze;
see Fig. 1.1.
The boundary and initial data are given as follows:

\smallskip
\begin{itemize}
\item[(\rmnum{1})] The Lipschitz function $g\in {\rm Lip}(\mathbb{R}_+; \mathbb{R})$ satisfies
$$
g(0) = g^{\prime}(0+) = 0, \hspace{5mm} \displaystyle\lim_{x\to\infty} \text{arctan}(g^{\prime}(x+)) = 0,
\hspace{5mm} g^{\prime} \in BV(\mathbb{R}_+; \mathbb{R}),
$$
and
$$
{\rm TV}(g^{\prime}(\cdot)) < \varepsilon \quad \hspace{5mm} \text{for some constant } \varepsilon > 0.
$$
Denote $\Omega \mathrel{\mathop:}= \lbrace (x, y)\,:\,y > g(x), x \geq 0 \rbrace$,
$\Gamma \mathrel{\mathop:}= \lbrace (x, y)\,:\,y = g(x), x \geq 0 \rbrace$,
and $\textbf{n}(x\pm)$ = $\frac{(-g^{\prime}(x\pm), 1)}{\sqrt{(g^{\prime}(x\pm))^2 + 1}}$ as
the outer normal vectors to $\Gamma$ at the respective points $x\pm$ ({\it cf.} Fig. \ref{Figure1}).

\item[(\rmnum{2})] The incoming flow $\overline{U}(y) := U_{0}^{b}(y) + \widetilde{U_0}(y)$
at $x=0$ is composed of two parts:

\smallskip
\begin{enumerate}
\item[(a)] The upstream flow $U_{0}^{b}(y)$ consists of one straight vortex sheet/entropy wave $y = y_{0}^{\ast} > 0$
and two constant vectors $U^b_0(y) = U_-$ when $0 < y < y_{0}^{\ast}$ and $U^b_0(y) = U_+$ when $y > y_{0}^{\ast} > 0$ such that
$$
v_- = v_+ = 0, \hspace{5mm} u_{\pm} > c_{\pm} > 0,
$$
where $c_{\pm} = \sqrt{\gamma p_{\pm}/\rho_{\pm}}$ is the sonic speed of state $U_{\pm}$.

\item[(b)] The $BV$ perturbation $\widetilde{U_0}(y) = (\tilde{\uu}_0, \tilde{p}_0, \tilde{\rho}_0)(y)\in (L^1 \cap BV)(\mathbb{R}; \mathbb{R}^{4})$
at $x = 0$ with ${\rm TV}(\widetilde{U_0}) \ll 1$.
\end{enumerate}
\end{itemize}

Then we consider the following initial-boundary value problem for system \eqref{2.1}:
\begin{align}
&\textbf{Boundary Condition:} && \uu \cdot \textbf{n}|_{\Gamma} = 0, \label{1.9}\\
&\textbf{Cauchy Condition:}   && U|_{x=0} = \overline{U}(y)=U_0^b(y)+\widetilde{U_0}(y).\quad\,\,\qquad \label{1.8}
\end{align}

\noindent
\textbf{Definition 1.1} (\textit{Admissible entropy solutions}).
A $BV$ function $U=U(x,y)$ is said to be an entropy solution of the initial-boundary value problem \eqref{2.1}--\eqref{1.8}
if and only if the following conditions hold:
\begin{itemize}
\item[(\rmnum{1})] $U$ is a weak solution of \eqref{2.1}
and satisfies \eqref{1.9}--\eqref{1.8}
in the trace sense;

\item[(\rmnum{2})] $U$ satisfies the {\it steady entropy Clausius inequality}:
\begin{equation}\label{1.10}
(\rho uS)_x + (\rho vS)_y \geq 0
\end{equation}
in the distributional sense in $\Omega$ including the Lipschitz wall boundary.
\end{itemize}

To solve the initial-boundary value problem \eqref{2.1}--\eqref{1.8}, in this paper, we develop suitable methods
to deal with the difficulties caused by the nonstrict hyperbolicity of the system
and the Lipschitz wall boundary, in comparison with the previous results for the Cauchy
problem of strictly hyperbolic systems of conservation laws.
For supersonic Euler flow with a strong shock-front emanating from the wedge vertex,
Chen-Li \ci{Chen-Li-2008} worked out the issue for a Lipschitz wedge boundary.
We now discuss some main differences in our work here from the Cauchy problem and the resulting
key difficulties. We remark that, in the case of the Cauchy problem concerning only $weak$ waves,
the decrease of the Lyapunov functional and the $L^1$--stability of the solutions were obtained
through the cancellation of distances on both sides of the waves.
In the presence of a strong shock, for the $L^1$--stability of solutions of the Cauchy problem
for strictly hyperbolic systems of conservation laws,
the Lyapunov functional was identified to decrease by employing the strength of the strong shock
to control the strengths of weak waves of the other families ({\it e.g.}, see Lewicka-Trivisa \ci{Lewicka-Trivisa-2002}).
In contrast with our Lipschitz wall problem, which is an initial-boundary value problem,
there is no such cancellation
as only one-side is possible near the boundary.
Furthermore, no strong vortex sheets/entropy waves (characteristic discontinuities) nor strong shocks
are present to handle the strength of the weak waves of the other families,
and the terms in the estimates for the first and fourth families carry different signs.
As such, it is difficult to say whether the functional can be made to decrease for our case of strong vortex
sheets and entropy waves with multiplicity of eigenvalues.
One of the key steps resolving this difficulty is to use the physical feature of the boundary
condition that the flow of two solutions near the boundary must run in parallel (also see \ci{Chen-Li-2008}).
This observation helps us obtain additional quantitative relations near the boundary.
Then, applying suitable weights and adjustments in the coefficients of the Lyapunov functional and using
the cancellation
between the different families, the functional is proved to decrease in the flow direction.

The rest of the paper is organized as follows:
In \S 2, we recall some fundamental properties of the $2$-D steady Euler system \eqref{1.1}, {\it i.e.}, \eqref{2.1},
and discuss related nonlinear waves and wave interaction estimates.
In \S 3, the wave-front tracking algorithm in the presence of strong vortex
sheets/entropy waves is discussed, the suitable interaction
potential $\mathcal{Q}$ including the effect of the Lipschitz wall is constructed,
and the existence of entropy solutions in $BV$ is established for the initial-boundary value problem.
In \S 4, we construct a
Lyapunov
functional
(equivalent to the $L^1$--distance between
two entropy solutions $U$ and $V$) to include the nonlinear waves produced by the wall boundary vertices.
Then, in \S 5, the monotone decrease of the Lyapunov functional
is established
in the flow direction, leading to the $L^1$--stability of the solutions containing
the strong vortex sheets/entropy waves. In \S 6, we employ the the estimates established in \S 3--\S 5
to obtain the existence of a Lipschitz semigroup
of solutions generated by a wave-front tracking approximation, as well as some estimates
on the uniformly Lipschitz semigroup
produced by the limit of wave-front tracking approximations.
Moreover,
the uniqueness of solutions with strong vortex sheets/entropy waves is obtained
in the larger class of viscosity solutions in \S 7.

\section{\small Steady Full Euler Equations: Nonlinear Waves and Wave Interactions}

In this section, we first present some basic properties of the steady Euler system \eqref{1.1}, {\it i.e.}, \eqref{2.1},
and then discuss nonlinear waves and related interaction estimates,
which will be employed in the subsequent development.

Notice that, when $U(x,y)\in C^1$,
system \eqref{2.1} is equivalent to
\begin{equation}
\nabla_{U}W(U)U_{x} + \nabla_{U}H(U)U_{y} = 0.
\end{equation}
Then the roots of the fourth degree polynomial
\begin{equation}
{\rm det}(\lambda\nabla_{U}W(U) - \nabla_{U}H(U))
\end{equation}
are the eigenvalues of \eqref{2.1}; that is, the solutions of the equation:
\begin{equation}
(v - \lambda u)^2\big((v - \lambda u)^2 - c^2(1 + \lambda^2)\big) = 0,
\end{equation}
where $c = \sqrt{\gamma p/\rho}$ is the sonic speed.
For supersonic flows ({\it i.e.}, $|\uu|> c$),
system \eqref{2.1} is hyperbolic.
Specifically, when $u > c$, system \eqref{2.1} has four real
eigenvalues in the $x$-direction:
\begin{align*}
\lambda_{j} &= \frac{uv + (-1)^{j} c\sqrt{u^2 + v^2 - c^2}}{u^2 - c^2} \qquad\mbox{for $j= 1,4$}; \nonumber \\
\lambda_{k} &= \frac{v}{u} \qquad \mbox{for $k =2,3$},
\end{align*}
with the four corresponding linearly independent eigenvectors given by
\begin{align}
\textbf{r}_{j} &= \kappa_{j}(-\lambda_{j}, 1, \rho(\lambda_{j}u - v),
\frac{\rho(\lambda_{j}u - v)}{c^2})^\top \qquad \mbox{for $j = 1,4$},  \label{2.3b} \\
\textbf{r}_{2} &= (u, v, 0, 0)^\top, \hspace{7mm} \textbf{r}_{3} = (0, 0, 0, \rho)^\top,\label{2.3a}
\end{align}
where $\kappa_{j}$ the re-normalization factors such
that $\textbf{r}_j \cdot \nabla \lambda_j =1$,
given that the $j$th-characteristic fields are genuinely nonlinear, $j = 1,4$.
The second and third linearly degenerate characteristic fields satisfy
$\textbf{r}_{k} \cdot \nabla \lambda_{k} = 0$, $k =2, 3$,
which correspond to vortex sheets and entropy waves, respectively.

\smallskip
The wave curves in the phase space for \eqref{2.1} are determined by
the Rankine-Hugoniot jump conditions:
\begin{equation}\label{2.5}
s[W(U)] =[H(U)],
\end{equation}
where $s$ is the propagation speed of the discontinuity.

There are two different waves associated with
the linearly degenerate families $\lambda_k=\frac{v_0}{u_0}, k=2,3$, with the corresponding
linearly independent right eigenvectors $\textbf{r}_k, k=2,3$,
in \eqref{2.3a}.

{\it Vortex sheets}:
\begin{equation}\label{2.6}
C_2(U_{0}):\hspace{3mm}  s=\frac{v}{u} =\frac{v_{0}}{u_{0}},
\hspace{4mm} p = p_{0}, \hspace{3mm} S=S_0, \hspace{3mm} u^2+v^2\ne u_0^2+v^2_0.
\end{equation}

{\it Entropy waves}:
\begin{equation}\label{2.7}
C_3(U_{0}):\hspace{3mm}  s=\frac{v}{u} =\frac{v_{0}}{u_{0}}, \hspace{4mm}
p = p_{0}, \hspace{3mm} (u,v)=(u_0,v_0), \hspace{3mm} S\ne S_0.
\end{equation}
The vortex sheet and the entropy wave
above match as a single characteristic discontinuity in the physical $(x,y)$--plane,
two independent parameters are needed to describe them in the phase
space $U = (\uu, p, \rho)=(u, v, p, \rho)$
since there are two linearly independent
eigenvectors corresponding to the repeated eigenvalues
$\lambda_{2} = \lambda_{3} =\frac{v}{u}$ of the linearly
degenerate characteristics fields.

The nonlinear waves associated with $\lambda_j, j=1,4$, are shock waves
and rarefaction waves. The propagation speeds of the shock waves are
$$
s_j\mathrel{\mathop:}= \frac{u_{0}v_{0} + (-1)^{j}\bar{c}_{0}\sqrt{u_{0}^2
+ v_{0}^2 - \bar{c}_{0}^2}}{u^2_0 - \bar{c}_{0}^2} \qquad \mbox{for $j =1,4$},
$$
where $\bar{c}_{0}^2 = \frac{c^{2}_{0}}{b_0} \frac{\rho}{\rho_{0}}$
and $b_{0} = \frac{\gamma + 1}{2} - \frac{\gamma - 1}{2} \frac{\rho}{\rho_{0}}$.
Substituting $s_j$ into \eqref{2.5},
the $j$-Hugoniot curves $S_j(U_{0})$ through state $U_{0}$ are
$$
S_j(U_0): \hspace{3mm} [p] = \frac{c^{2}_{0}}{b_{0}}[\rho],
\hspace{3mm} [u] = -s_j[v],
\hspace{3mm} \rho_{0}(s_ju_{0} - v_{0})[v] = [p] \qquad \mbox{for $j = 1,4$}.
$$
Written as $S_j^{+}(U_{0})$, $j= 1, 4$,
the half curves of $S_j(U_{0})$ for $\rho > \rho_{0}$ in the phase space are said
to be the shock curves on which any state forms a shock with the below state $U_{0}$
in the $(x,y)$--plane respecting the entropy condition \eqref{1.10}.
Furthermore, for each $j=1$ or $j=4$,
curves $S^{+}_j(U_{0})$ and $R^{-}_j(U_{0})$ at state $U_{0}$
have the same curvature.

If $U$ is a piecewise smooth solution,
then any of the following conditions below is equivalent to
the entropy inequality \eqref{1.10} in Definition 1.1 for a shock wave (see also \ci{Chen-Zhang-Zhu-2007}):
\begin{itemize}
\item[(\rmnum{1})] {\it The physical entropy condition}: The density increases across the shock in the flow direction,
\begin{equation}\label{2.8}
 \rho_{\text{back}} < \rho_{\text{front}}.
\end{equation}
\item[(\rmnum{2})] {\it The Lax entropy condition}: On the $j$th-shock, the shock speed $s_j$ satisfies
\begin{eqnarray}\label{2.9-2.10}
&&\lambda_j(\text{back}) < s_j < \lambda_j(\text{front}) \hspace{5mm}\mbox{for $j = 1,4$}, \\
&& s_1 < \lambda_{2,3}(\text{back}), \hspace{5mm} \lambda_{2,3}(\text{front}) < s_{4}.
\end{eqnarray}
\end{itemize}

The rarefaction wave curves $R_{j}^{-}(U_{0})$ through state $U_{0}$ in the state space are given by
\begin{equation}\label{2.11}
R^{-}_j: \hspace{3mm} {\rm d}p = c^2 {\rm d}\rho, \hspace{3mm} {\rm d}u = -\lambda_j{\rm d}v, \hspace{3mm}
\rho(\lambda_ju -v){\rm d}v = {\rm d}p \hspace{8mm} \text{ for $\rho < \rho_{0}$}, \quad j = 1, 4.
\end{equation}

We now discuss several essential properties of the nonlinear waves and related wave interaction estimates
in Lemmas 2.1--2.7 below.
These facts will be used in the subsequent development;
see also Chen-Zhang-Zhu \ci{Chen-Zhang-Zhu-2007}
for further details.
\vspace{3mm}

\noindent
\emph{2.1. Riemann Problems and Riemann Solutions}

\smallskip
We focus on the related Riemann problems and their solutions in this section,
which serve as the building blocks for the front tracking algorithm
for the initial-boundary value problem \eqref{2.1}--\eqref{1.8}.

\vspace{2.5mm}
\noindent
\emph{Lateral Riemann problem.}
Consider the following lateral Riemann problem with boundary $\Gamma$:
\begin{equation}\label{2.14a}
\begin{cases}
(1.4), \hspace{2.2cm} &{} \\
U|_{x=x_0} = \underline{U}, & \\
\uu \cdot \textbf{n}|_{\Gamma} = 0.
\end{cases}
\end{equation}
It has been observed in \ci{Courant-Friedrichs-1948} that, if the angle between the flow direction of
the constant front-state $\underline{U}$ and the wall at a boundary vertex is smaller than $\pi$
and larger than the extreme angle determined by the incoming flow state and $\gamma > 1$,
then a unique $4$-shock is generated, separating the front-state from the supersonic back-state.
If the angle between the flow direction of the front-state and the wall at a boundary vertex
is larger than $\pi$ and less than the extreme angle,
then a $4$-rarefaction wave is produced, emanating from the vertex.
These waves are easily seen through the shock polar analysis
({\it cf.} \ci{Chen-Zhang-Zhu-2007, Courant-Friedrichs-1948}).
This signifies that, when the angle between the flow direction of the front-state and the wall
at a boundary vertex is close to $\pi$, the lateral Riemann problem can be uniquely solved.
For further details, see
\ci{Chen-Zhang-Zhu-2007,Courant-Friedrichs-1948}.

In particular, the background solution $U=U_0^b$ is the unique entropy solution of
problem \eqref{2.14a} with $\underline{U}=U_0^b$ and $g \equiv 0$, consisting of
two constant states $U_{-}= (u_{-},0,p_{-},\rho_{-})$
and $U_{+} = (u_{+},0,p_{+},\rho_{+})$, satisfying $u_{\pm} > c_{\pm} > 0$ in
subdomains $\Omega_+$ and $\Omega_-$ of $\Omega$ separated by the straight vortex sheet/entropy wave.
The principal aim of this paper is to establish the $L^1$ well-posedness
for problem \eqref{2.1}--\eqref{1.8}
for the solutions near the background solution $U_0^b$ containing the strong vortex sheet/entropy wave
$\lbrace U_{-},U_{+}\rbrace$.

\vspace{3mm}
\noindent
\emph{Riemann problem involving only weak waves.}
Consider the Riemann problem:
\begin{equation}\label{2.12}
\begin{cases} (1.4),
\\[2mm]
U|_{x=x_{0}} = \underline{U} =
\begin{cases}
U_{a} \hspace{2mm}\mbox{for $y > y_{0}$}, \\
U_{b} \hspace{2mm}\mbox{for $y < y_{0}$},
\end{cases}
\end{cases}
\end{equation}
with the constant states $U_{a}$ and $U_{b}$ denoting the {\it above} state
and {\it below} state with respect to line $y = y_{0}$, respectively.
Then there exists $\varepsilon > 0$ so that, for any states $U_{b}$ and $U_{a}$ in the neighborhood
$\textit{O}_{\varepsilon}(U_{+})$ of $U_{+}$,
or $U_{b}$ and $U_{a}$ in the neighborhood $\textit{O}_{\varepsilon}(U_{-})$ of $U_{-}$,
the Riemann problem \eqref{2.12} has a unique admissible solution
consisting of at most four waves of shocks, rarefaction waves, a vortex sheet, and an entropy wave.

\vspace{3mm}
\noindent
\emph{Riemann problem involving the strong vortex sheets/entropy waves.}
From now on, the notation $\lbrace U_{b}, U_{a}\rbrace = \left(\alpha_1, \alpha_2, \alpha_3, \alpha_4 \right)$
will be used to write $U_{a} = {\Phi}\left( \alpha_4, \alpha_3, \alpha_2, \alpha_1; U_{b}\right)$
as the solution of the Riemann problem,
where ${\Phi} \in C^2$, and $\alpha_{j}$ is the strength of the $j$-wave (measuring the jump across the wave).
For any waves with $U_{b} \in O_{\varepsilon}(U_{-})$ and $U_{a} \in O_{\varepsilon}(U_{+})$,
we also use $\lbrace U_{b}, U_{a} \rbrace = \left(0, \sigma_{2}, \sigma_{3}, 0\right)$
to denote the strong vortex sheet/entropy wave that connects $U_{b}$ and $U_{a}$ with
strength $\left(\sigma_{2}, \sigma_{3}\right)$. That is,
$$
U_{m} = \Phi_2(\sigma_{2}; U_{b})
   \mathrel{\mathop:}= (u_{b}e^{\sigma_{2}}, v_{b}e^{\sigma_{2}}, p_{b}, \rho_{b}),
\qquad
U_{a} = \Phi_3(\sigma_{3}; U_{m}) \mathrel{\mathop:}
= (u_{m}, v_{m}, p_{m}, \rho_{m}e^{\sigma_{3}}).
$$
In particular, for the background solution $U_0^b$, $\{U_-, U_+\}=(0, \sigma_{20}, \sigma_{30}, 0)$:
$$
U_{+} = (u_{+}, 0, p_{+}, \rho_{+}) = (u_{-}e^{\sigma_{20}}, 0, p_{-}, \rho_{-}e^{\sigma_{30}}).
$$
We write $G(\sigma_{3}, \sigma_{2}; U_{b}):= \Phi_3(\sigma_{3}; \Phi_2(\sigma_{2}; U_{b}))$
for any $U_{b} \in O_{\varepsilon}(U_{-})$. Then we have

\medskip
\noindent
$\textbf{Lemma 2.1.}$
\emph{The vector function $G(\sigma_{3}, \sigma_{2}; U_{b})$ satisfies}
\begin{equation}\label{2.13}
G_{\sigma_{2}}\left(\sigma_{3}, \sigma_{2}; U_{b}\right)
= \left(u_{b}e^{\sigma_{2}}, v_{b}e^{\sigma_{2}}, 0, 0\right),
\hspace{4mm} G_{\sigma_{3}}\left(\sigma_{3}, \sigma_{2}; U_{b}\right) = (0, 0, 0, \rho_{b}e^{\sigma_{3}}),
\end{equation}
\emph{and}
\begin{equation}\label{2.14}
\nabla_{U}G(\sigma_{3}, \sigma_{2}; U_{b}) = {\rm diag}(e^{\sigma_{2}}, e^{\sigma_{2}}, 1, e^{\sigma_{3}}).
\end{equation}
\emph{Furthermore, for the background plane vortex sheet and entropy wave with the below
state $U_{-} = (u_{-}, 0, p_{-}, \rho_{-})$,
above state $U_{+} = (u_{+}, 0, p_{+}, \rho_{+})$, and strength $(\sigma_{20}, \sigma_{30})$},
\begin{equation}\label{2.15}
{\rm det}(\textbf{r}_4(U_{+}), G_{\sigma_{3}}(\sigma_{30}, \sigma_{20}; U_{-}),
G_{\sigma_{2}}(\sigma_{30}, \sigma_{20}; U_{-}), \nabla_{U}G(\sigma_{30}, \sigma_{20};
U_{-}) \cdot \textbf{r}_{1}(U_{-})) > 0.
\end{equation}

These can be easily obtained from direct calculations, which are thus omitted.

The properties in \eqref{2.13}--\eqref{2.15} above play a fundamental role
in achieving the necessary estimates for the strengths of reflected weak waves
in the interaction between the strong vortex sheet/entropy wave
and weak waves (see the proofs for Lemmas 2.4--2.7).
\vspace{4mm}

\noindent
\emph{2.2. Estimates for Wave Interactions and Reflections.}
In the following, several essential estimates for wave interactions and reflections are provided.
For their proofs
and related details,
see \ci{Chen-Zhang-Zhu-2007}.

\smallskip
{\it Estimates for weak wave interactions.} For the weak wave interaction away from both
the strong vortex sheet/entropy wave
and the wall boundary in subdomains $\Omega_+$ or $\Omega_-$, we have the following estimate:

\smallskip
\noindent
$\textbf{Lemma 2.2.}$
\emph{Assume that
$U_{b}, U_{m}, U_{a} \in O_{\varepsilon}(U_{+})$, or $U_{b}, U_{m}, U_{a} \in O_{\varepsilon}(U_{-})$,
are three states with $\lbrace U_{b}, U_{m}\rbrace = (\alpha_{1}, \alpha_{2}, \alpha_{3}, \alpha_{4})$
and $\lbrace U_{m}, U_{a}\rbrace = (\beta_{1}, \beta_{2}, \beta_{3}, \beta_{4})$.
Then $\lbrace U_{b}, U_{a}\rbrace = (\gamma_{1}, \gamma_{2}, \gamma_{3}, \gamma_{4})$ with}
\begin{equation}\label{lemma2.2}
\gamma_i = \alpha_{i} + \beta_{i} + O(1)\Delta(\alpha, \beta),
\end{equation}
\emph{where $\Delta(\alpha, \beta) = (|\alpha_{4}| + |\alpha_{3}| + |\alpha_{2}|)|\beta_{1}|
+ |\alpha_{4}|(|\beta_{3}| + |\beta_{2}|) + \sum_{j=1,4}\Delta_{j}(\alpha, \beta)$ and}
\begin{equation}\label{2.16}
\Delta_{j}(\alpha, \beta) =
\begin{cases}
0, \hspace{3mm} &\alpha_{j} \geq 0, \beta_{j} \geq 0, \\
|\alpha_{j}||\beta_{j}|, \hspace{2mm}& otherwise.
\end{cases}
\end{equation}

\smallskip
{\it Estimates on the boundary perturbation of weak waves and the reflection of weak waves on the boundary.}
We write $\lbrace C_{l}(a_{l}, b_{l})\rbrace^{\infty}_{l=0}$ for
points $\lbrace(a_{l},b_{l})\rbrace^{\infty}_{l=0}$ in the $(x,y)$--plane with $0 < a_{l} < a_{l+1}$.
Define
\begin{equation}\label{2.17}
\begin{cases}
\theta_{{l},{l+1}} = \text{arctan}\left(\frac{b_{l+1} - b_{l}}{a_{l+1} - a_{l}}\right),
\hspace{3mm} \theta_{l} = \theta_{{l},{l+1}} - \theta_{{l-1},{l}}, \hspace{3mm} \theta_{-1, 0} =0,   \\[2mm]
\Omega_{l+1} = \lbrace (x, y)\,:\, x \in \left[a_{l}, a_{l+1}\right],
y > b_{l} + (x - a_{l})\text{tan}(\theta_{{l}, {l+1}})\rbrace, \\[2mm]
\Gamma_{l+1} = \lbrace (x, y)\,:\, x \in \left(a_{l}, a_{l+1}\right),
y = b_{l} + (x - a_{l})\text{tan}(\theta_{{l}, {l+1}})\rbrace,
\end{cases}
\end{equation}
and the outer normal vector to $\Gamma_{l}$:
\begin{equation}\label{2.18}
\textbf{n}_{l+1} = \frac{(b_{l+1} - b_{l}, a_{l} - a_{l+1})}{\sqrt{(b_{l+1} - b_{l})^2 + (a_{l+1} - a_{l})^2}}
= (\text{sin}(\theta_{l,l+1}), -\text{cos}(\theta_{l,l+1})).
\end{equation}

With the constant state $\underline{U}$, consider the following lateral Riemann problem:
\begin{equation}\label{2.19}
\begin{cases}
(2.1) \hspace{2.2cm} &\text{in $\Omega_{l+1}$}, \\
U|_{x=a_{l}} = \underline{U}, \\
\uu \cdot \textbf{n}_{l+1}|_{\Gamma_{l+1}} = 0.
\end{cases}
\end{equation}

\noindent $\textbf{Lemma 2.3.}$
\emph{Suppose $\left\lbrace U_{m}, U_{a} \right\rbrace = (\beta_{1}, \beta_{2}, \beta_{3}, 0)$
and $\left\lbrace U_{l}, U_{m}\right\rbrace = (0, 0, 0, \alpha_{4})$ with
$\uu_{l}\cdot \textbf{n}_{l}|_{\Gamma_l} = 0$}.
\emph{Then there exists a unique solution $U_{l+1}$ of problem {\eqref{2.19}} such that
$\left\lbrace U_{l+1}, U_{a}\right\rbrace = (0, 0, 0, \delta_{4})$
and $\uu_{l+1} \cdot \textbf{n}_{l+1}|_{\Gamma_{l+1}}= 0$. Moreover,}
\begin{equation}\label{2.20}
\delta_{4} = \alpha_{4} + K_{b1}\beta_{1} + K_{b2}\beta_{2} + K_{b3}\beta_{3} + K_{b0}\theta_{l},
\end{equation}
\emph{where $K_{b1}$, $K_{b2}$, $K_{b3}$, and $K_{b0}$ are $C^{2}$--functions of $\beta_{3}$, $\beta_{2}$, $\beta_{1}$,
$\alpha_{4}$, $\theta_{l+1}$, and $U_{a}$ satisfying}
\begin{equation}\label{2.21}
K_{b1}|_{\left\lbrace \theta_{l} = \alpha_{4} = \beta_{1} = \beta_{2} = \beta_{3}=0,
 U_{a} = U_{-}\right\rbrace} = 1,\hspace{0.8cm}
 K_{bi}|_{\left\lbrace\theta_{l}= \alpha_{4}= \beta_{1}= \beta_{2}= \beta_{3}=0,
 U_{a} = U_{-}\right\rbrace} = 0 \,\,\,\,\mbox{for $i = 2, 3$},
\end{equation}
\emph{and $K_{b0}$ is bounded. In particular, $K_{b0} < 0$ at the origin.}

\smallskip
This lemma has two purposes. The first is to estimate the weak waves generated by the vertices
on the Lipschitz wall boundary.
This boundedness will be used to control the boundary perturbation; see \eqref{3.2} below in the construction
of the wave interaction potential $\mathcal{Q}(x)$.
The second is to estimate the strength of the reflected wave $\delta_{4}$ with respect to the incident
wave $\alpha_{1}$. Property \eqref{2.21} of the coefficients will play an important role in controlling
the reflected waves.
\vspace{3mm}

{\it Estimates for the interactions between the strong vortex sheet/entropy wave
and weak waves from below.}
Estimate \eqref{key-1} below plays a key role in ensuring the $L^{1}$--stability of entropy solutions,
especially for the existence of constants $w^{b}_{1}$ and $w^{b}_{4}$ in Lemma 5.1
(see below). This estimate also ensures the existence of $K^{\ast} \in (K_{11}, 1)$ in the
construction of the wave interaction potential $\mathcal{Q}(x)$ in \eqref{3.2}.

\smallskip
\noindent $\textbf{Lemma 2.4.}$ \emph{Let $U_{b}, U_{m} \in O_{\varepsilon} (U_{-})$
and $U_{a} \in O_{\varepsilon}(U_{+})$ with}
$$
\{U_{b}, U_{m}\} = (0, \alpha_{2}, \alpha_{3}, \alpha_{4}),
\hspace{7mm} \{U_{m}, U_{a}\} = (\beta_{1}, \sigma_{2}, \sigma_{3}, 0).
$$
\emph{Then there exists a unique $(\delta_{1}, \sigma^{\prime}_{2}, \sigma^{\prime}_{3}, \delta_{4})$ such that the Riemann problem \eqref{2.12} admits an admissible solution that consists of a weak $1-$wave
of strength $\delta_{1}$, a strong vortex sheet/entropy wave of strength $(\sigma^{\prime}_{2}, \sigma^{\prime}_{3})$, and a weak $4-$wave of strength $\delta_{4}$}:
$$
\{U_{b}, U_{a}\} = (\delta_{1}, \sigma^{\prime}_{2}, \sigma^{\prime}_{3}, \delta_{4}),
$$
so that
\begin{eqnarray}
&&\delta_{1} = \beta_{1} + K_{11}\alpha_{4} + O(1)\Delta^{\prime},
\hspace{5mm} \delta_{4} = K_{14}\alpha_{4} + O(1)\Delta^{\prime},\nonumber\\
&& \sigma_{2}^{\prime} = \sigma_{2} + \alpha_{2} + K_{12}\alpha_{4} + O(1)\Delta^{\prime},
\hspace{5mm} \sigma_{3}^{\prime} = \sigma_{3} + \alpha_{3} + K_{13}\alpha_{4} + O(1)\Delta^{\prime},\nonumber\\
&&|K_{11}|_{\left\lbrace\alpha_{4}=\alpha_{3}=\alpha_{2}=0,\sigma_{2}=\sigma_{20},
\sigma_{3}=\sigma_{30}\right\rbrace} = \left|\frac{\lambda_{4}
(U_{+})e^{2\sigma_{20} + \sigma_{30}} - \lambda_{4}(U_{-})}{\lambda_{4}(U_{+})e^{2\sigma_{20} +\sigma_{30}} + \lambda_{4}(U_{-})}\right| < 1,
\label{key-1}
\end{eqnarray}
\emph{where $\sum_{j=2}^{4}|K_{1j}|$ is bounded, and $\Delta^{\prime} = |\beta_{1}|(|\alpha_{2}| + |\alpha_{3}|)$.}
\vspace{4mm}

\noindent $\textbf{Lemma 2.5.}$ \emph{The coefficient,
$|K_{14}|_{\lbrace\alpha_{4}=\alpha_{3}=\alpha_{2}=0,\sigma_{2}=\sigma_{20},
\sigma_{3}=\sigma_{30}\rbrace}$,
in the strength $\delta_4$ of a weak $4$-wave in Lemma {\rm 2.4}
remains bounded away from zero.}
\vspace{4mm}

\noindent $\textbf{Proof.}\,$  By Lemma 2.4, we can find a unique solution
$(\delta_1, \sigma^{\prime}_2, \sigma^{\prime}_3, \delta_4)$ as a $C^2$--function
of $\alpha_2, \alpha_3, \alpha_4, \beta_1,$ $\sigma_2, \sigma_3$, and $U_b$ to
\begin{equation}\label{2.23}
\Phi_4(\delta_4; G(\sigma^{\prime}_3, \sigma^{\prime}_2; \Phi_1(\delta_1; U_b)))
= G(\sigma_3, \sigma_2; \Phi_1(\beta_1; \Phi(\alpha_4, \alpha_3, \alpha_2, 0; U_b))).
\end{equation}
That is,
$$
\sigma^{\prime}_i = \sigma^{\prime}_i(\alpha_2,\alpha_3, \alpha_4, \beta_1, \sigma_2, \sigma_3)
\hspace{3mm} \mbox{for $i=2,3$}, \hspace{9mm}
\delta_{j} = \delta_{j}(\alpha_2, \alpha_3, \alpha_4, \beta_1, \sigma_2, \sigma_3) \hspace{3mm} \mbox{for $j=1,4$},
$$
where we have omitted $U_b$ for simplicity.
Moreover, from \ci{Chen-Zhang-Zhu-2007}, we have
\begin{equation}
K_{1j} = \int\limits_{0}^{1}\ \partial_{\alpha_4}\delta_j(\alpha_2, \alpha_3,
\theta\alpha_4, \beta_1, \sigma_2, \sigma_3) \,{\rm d}\theta
\qquad\,\, \mbox{for $j= 1,4$}. \nonumber
\end{equation}

Differentiate \eqref{2.23} with respect to $\alpha_4$,
and let $\beta_1 = \alpha_4 = \alpha_3 = \alpha_2 = 0, \sigma_2 = \sigma_{20}$, and $\sigma_3 = \sigma_{30}$.
We obtain
\begin{align*}
\nabla_{U}G(\sigma_{30}, \sigma_{20}; U_-) \cdot \textbf{r}_4(U_-)
&=\partial_{\alpha_4}\delta_{4}\,\textbf{r}_{4}(U_+)
+ \partial_{\alpha_4}\sigma^{\prime}_{3}\,G_{\sigma_{3}}(\sigma_{30}, \sigma_{20}; U_-) \\
&\quad + \partial_{\alpha_4}\sigma^{\prime}_{2}\,G_{\sigma_{2}}(\sigma_{30}, \sigma_{20}; U_-)
  + \partial_{\alpha_4}\delta_{1}\,\nabla_{U}G(\sigma_{30}, \sigma_{20}; U_-) \cdot \textbf{r}_1(U_-).
\end{align*}

\noindent By Lemma 2.1, we have \\

\noindent
{\small
\begin{align*}
&|\partial_{\alpha_4}\delta_4|\\[1mm]
&= \left|\frac{ \text{det}(\nabla_{U}G(\sigma_{30}, \sigma_{20}; U_-) \cdot \textbf{r}_{4}(U_-),
G_{\sigma_{3}}(\sigma_{30}, \sigma_{20}; U_-), G_{\sigma_{2}}(\sigma_{30}, \sigma_{20}; U_-),
\nabla_{U}G(\sigma_{30}, \sigma_{20}; U_-) \cdot \textbf{r}_1(U_-))}{ \text{det}(\textbf{r}_{4}(U_+), G_{\sigma_{3}}(\sigma_{30}, \sigma_{20}; U_-), G_{\sigma_{2}}(\sigma_{30}, \sigma_{20}; U_-), \nabla_{U}G(\sigma_{30}, \sigma_{20}; U_-) \cdot \textbf{r}_{1}(U_-))}\right|\\[1mm]
&= \left|\frac{\kappa_1(U_-)\kappa_4(U_-)\rho^{2}_{-}u^{2}_{-}e^{2\sigma_{20}+\sigma_{30}}(\lambda_{4}(U_{-}) - \lambda_{1}(U_{-}))}{\kappa_1(U_-)\kappa_4(U_+)\rho^{2}_{-}u^{2}_{-}e^{\sigma_{20}
+\sigma_{30}}(\lambda_{4}(U_{+})e^{2\sigma_{20}+\sigma_{30}} - \lambda_{1}(U_{-}))}\right| \\[1mm]
&= \left| \frac{2\kappa_{4}(U_{-})e^{\sigma_{20}}\lambda_4(U_-)}{\kappa_{4}(U_+)(\lambda_{4}(U_+)e^{2\sigma_{20}+\sigma_{30}} + \lambda_{4}(U_{-}))}\right| > 0.
\end{align*}
}
This completes the proof.

\vspace{4mm}
\emph{Estimates for the interactions between the strong vortex sheet/entropy wave and weak waves from above.} We have

\smallskip
\noindent
$\textbf{Lemma 2.6.}$
\emph{Let $U_{b} \in O_{\varepsilon}(U_{-})$ and $U_{m}, U_{a} \in O_{\varepsilon}(U_{+})$ with}
$$
\lbrace U_{b}, U_{m}\rbrace = (0, \sigma_{2}, \sigma_{3}, \alpha_{4}),
\hspace{7mm} \lbrace U_{m}, U_{a}\rbrace = (\beta_{1}, \beta_{2}, \beta_{3}, 0).
$$
\emph{Then there exists a unique $(\delta_{1}, \sigma^{\prime}_{2}, \sigma^{\prime}_{3}, \delta_{4})$
such that the Riemann problem \eqref{2.12} admits an admissible solution that consists of a weak $1-$wave
of strength $\delta_{1}$, a strong vortex sheet/entropy wave of strength $(\sigma^{\prime}_{2}, \sigma^{\prime}_{3})$,
and a weak $4-$wave of strength $\delta_{4}$}:
$$
\lbrace U_{b}, U_{a}\rbrace = (\delta_{1}, \sigma^{\prime}_{2}, \sigma^{\prime}_{3}, \delta_{4}),
$$
so that
\begin{align*}
&\delta_{1} = K_{21}\beta_{1} + O(1)\Delta^{''}, \hspace{5mm}
\sigma_{2}^{\prime} = \sigma_{2} + \beta_{2} + K_{22}\beta_{1} + O(1)\Delta^{''},\\
&\sigma_{3}^{\prime} = \sigma_{3} + \beta_{3} + K_{23}\beta_{1} + O(1)\Delta^{''},
\hspace{5mm} \delta_{4} = \alpha_{4} + K_{24}\beta_{1} + O(1)\Delta^{''},
\end{align*}
\emph{where $\sum_{j=1}^{4}|K_{2j}|$ is bounded, and $\Delta^{''} = |\alpha_{4}|(|\beta_{2}| + |\beta_{3}|).$}
\vspace{4mm}

The constant, $K_{21}$, here is used in the definition of weighted strength $b_{\alpha}$ of weak waves in \eqref{3.1}.

\vspace{4mm}
\noindent
$\textbf{Lemma 2.7.}$ \emph{The coefficient,
$|K_{21}|_{\lbrace\beta_{3}=\beta_{2}=\beta_{1}=0,
\sigma_{2}=\sigma_{20},\sigma_{3}=\sigma_{30}\rbrace}$,
in the strength $\delta_1$ of a weak $1$-wave in Lemma {\rm 2.6} remains bounded away from zero,
while the reflection coefficient $|K_{24}|<1$.}
\vspace{4mm}

\noindent $\textbf{Proof.}\,$
By Lemma 2.6, we can find a unique solution $(\delta_1, \sigma^{\prime}_2, \sigma^{\prime}_3, \delta_4)$
as a $C^2$--function of $\alpha_2, \alpha_3, \alpha_4$, $\beta_1$, $\sigma_2$, $\sigma_3$, and $U_b$ to
\begin{equation}\label{2.24}
\Phi_4(\delta_4; G(\sigma^{\prime}_3, \sigma^{\prime}_2; \Phi_1(\delta_1; U_b)))
= \Phi(0, \beta_3, \beta_2, \beta_1; \Phi_4(\alpha_4; G(\sigma_3, \sigma_2; U_b))).
\end{equation}
That is,
$$
\sigma^{\prime}_i = \sigma^{\prime}_i(\beta_1,\beta_2, \beta_3, \alpha_4, \sigma_2, \sigma_3)
 \hspace{3mm} \mbox{for $i=2,3$},
\hspace{9mm} \delta_{j} = \delta_{j}(\beta_1,\beta_2, \beta_3, \alpha_4, \sigma_2, \sigma_3)
 \hspace{3mm} \mbox{for $j=1,4$},
$$
where we have omitted $U_b$ for simplicity.
Moreover, from \ci{Chen-Zhang-Zhu-2007}, we have
\begin{equation}
K_{2j} = \int\limits_{0}^{1}\ \partial_{\beta_1}\partial_j(\theta\beta_1, \beta_2, \beta_3,
\alpha_4, \sigma_2, \sigma_3) \,{\rm d}\theta \hspace{7mm} \mbox{for $j= 1,4$}. \nonumber
\end{equation}

Similarly, differentiate \eqref{2.24} with respect to $\beta_1$, and
let $\alpha_4 = \beta_1 = \beta_2 = \beta_3 = 0, \sigma_2 = \sigma_{20}$, and $\sigma_3 = \sigma_{30}$.
Then we obtain
\begin{align*}
\textbf{r}_1(U_+)
&=\partial_{\beta_1}\delta_{4}\,\textbf{r}_{4}(U_+)
+ \partial_{\beta_1}\sigma^{\prime}_{3}\,G_{\sigma_{3}}(\sigma_{30}, \sigma_{20}; U_-) \\
&\quad+ \hspace{1mm} \partial_{\beta_1}\sigma^{\prime}_{2}\,G_{\sigma_{2}}(\sigma_{30}, \sigma_{20}; U_-) + \partial_{\beta_1}\delta_{1}\,\nabla_{U}G(\sigma_{30}, \sigma_{20}; U_-) \cdot \textbf{r}_1(U_-).
\end{align*}

\noindent By Lemma 2.1, we have
{\small
\begin{align*}
|\partial_{\beta_1}\delta_1|
&=\left|\frac{\text{det}(\textbf{r}_{4}(U_+), G_{\sigma_{3}}(\sigma_{30}, \sigma_{20}; U_-),
G_{\sigma_{2}}(\sigma_{30}, \sigma_{20}; U_-),\textbf{r}_1(U_+))}{ \text{det}(\textbf{r}_{4}(U_+), G_{\sigma_{3}}(\sigma_{30}, \sigma_{20}; U_-), G_{\sigma_{2}}(\sigma_{30}, \sigma_{20}; U_-), \nabla_{U}G(\sigma_{30}, \sigma_{20}; U_-) \cdot \textbf{r}_{1}(U_-))}\right|\\[1mm]
&= \left| \frac{2\kappa_{1}(U_{+})\lambda_4(U_+)e^{\sigma_{20}+\sigma_{30}}}{\kappa_{1}(U_-)(\lambda_{4}(U_+)e^{2\sigma_{20}+\sigma_{30}}
+ \lambda_{4}(U_{-}))}\right| > 0.
\end{align*}
}

\noindent However, for the reflection coefficient $|K_{24}|$, we have
\vspace{4mm}
\begin{align*}
|\partial_{\beta_1}\delta_4|
&= \left|\frac{\text{det}(\textbf{r}_1(U_+), G_{\sigma_{3}}(\sigma_{30}, \sigma_{20}; U_-),
G_{\sigma_{2}}(\sigma_{30}, \sigma_{20}; U_-), \nabla_{U}G(\sigma_{30}, \sigma_{20}; U_-) \cdot \textbf{r}_1(U_-))}{ \text{det}(\textbf{r}_{4}(U_+), G_{\sigma_{3}}(\sigma_{30}, \sigma_{20}; U_-), G_{\sigma_{2}}(\sigma_{30}, \sigma_{20}; U_-), \nabla_{U}G(\sigma_{30}, \sigma_{20}; U_-) \cdot \textbf{r}_{1}(U_-))}\right|  \\[1mm]
&= \left| \frac{-\lambda_4(U_+)e^{2\sigma_{20}+\sigma_{30}}+\lambda_4(U_-)}{\lambda_{4}(U_+)e^{2\sigma_{20}+\sigma_{30}} + \lambda_{4}(U_{-})}\right| < 1,
\end{align*}
where $|K_{24}|$ is not necessarily bounded away from zero, but is less than one.

\smallskip
\section{\small The Wave-Front Tracking Algorithm and Global Existence of Entropy Solutions}

We first start with a brief description of the wave-front tracking method
to be employed throughout in \S 4--\S 7, and then establish the existence of entropy solutions
when the perturbation of the incoming flow has small total variation at $x = 0$.

The main scheme in the wave-front tracking method is to construct approximate solutions
within a class of piecewise constant functions.
We first approximate the initial data by a piecewise constant vector function.
Then we solve the resulting Riemann problems exactly, with the exception of the rarefaction waves
that are replaced by the rarefaction fans with many small wave-fronts of equal strengths. The outgoing
fronts are continued up to the first time when two waves collide and a new Riemann problem is solved.
In this process, one has to modify the algorithm and introduce a simplified Riemann solver
in order to keep the number of wave-fronts finite for all $x \geq 0$ in the flow direction.
See Bressan \cite{Bressan-1992,Bressan-2000}
and Baiti-Jenssen \cite{Baiti-Jenssen-1998} for related references.

\vspace{3mm}
\noindent
\emph{3.1. The Riemann Solvers}.

\vspace{2mm}
As indicated in \S 2, the solution of the Riemann problem $\lbrace U_b, U_a \rbrace$
is a self-similar solution given by at most five states separated by shocks, vortex sheet/entropy wave,
or rarefaction waves.
To connect state $U_a$ to $U_b$, there exist $C^2$--curves $\eta \rightarrow \varphi(\eta)(U)$
with parametrization (which is equivalent to the arc length and consistent with renormalization $\mathbf{r}_j\cdot \nabla\lambda_j=1, j=1,4$)
such that
$$
U_b =\varphi(\eta)(U_a):= \Upsilon_4(\eta_4) \circ \cdots \circ \Upsilon_1(\eta_1)(U_a)
$$
for some $\eta = (\eta_1, \ldots, \eta_4)$,
and $U_j = \Upsilon_j(\eta_j) \circ \cdots \circ \Upsilon_1(\eta_1)(U_a)$ for $j = 1,2,3$.

\medskip
Next, we describe the construction of front tracking approximations for
the initial-boundary value problem \eqref{2.1}--\eqref{1.8}.
Denote $\vartheta>0$ as the initial approximation parameter.
Then the given initial data function $\overline{U}$
is first approximated by
a sequence of piecewise constant functions $\overline{U}^{\vartheta}$
in the $L^{1}$--norm,
and the wall boundary is also approximated as described in \eqref{2.17} in \S 2
with
$$
a_l=l\Delta x, \quad b_l=g(l\Delta x)  \qquad\,\, \mbox{for some $\Delta x>0$}.
$$
For fixed $\vartheta>0$, denote $\mathcal{Z}_{\vartheta}$ as the set of the total number of jumps in the
approximate initial data functions $\overline{U}^{\vartheta}$
and the tangential angle function of the wall boundary.
Let $\delta_{\vartheta} > 0$ be a parameter so that a rarafaction wave is replaced by a step function
whose {\it steps} are no further apart than $\delta_{\vartheta}$.
The discontinuity between two steps is set to propagate with a speed equal to the Rankine-Hugoniot speed
of the jump connecting the states corresponding to the two steps.
At any time, the simplified Riemann solver (defined below) is employed
with constant $\hat{\lambda}>0$ (as the speed of the generated non-physical wave)
which is strictly greater than all the wave speeds of system (2.1).
The strength of the non-physical wave is the error generated when
the simplified Riemann solver is applied.

\smallskip
$\textbf{Accurate Riemann solver.}$ The accurate Riemann solver ({\bf ARS}) is the exact solution
of the Riemann problem,
with the exception that every rarefaction wave $\lbrace w, R_j(w)(\alpha)\rbrace$, $j = 1, 4$,
is divided into equal parts and replaced by a piecewise constant rarefaction fan of several
new wave-fronts of equal strength.

\smallskip
$\textbf{Simplified Riemann solver.}$  When only very weak waves are involved,
the simplified Riemann solver ($\textbf{SRS}$) here is the same as the one described
in \ci{Baiti-Jenssen-1998, Bressan-2000}.
That is, all new weak waves are put together in a single $\textit{non-physical front}$ with
positive speed larger than all the characteristic speeds.
In the case of a weak wave interacting with the strong vortex sheet/entropy wave,
the purpose of ($\textbf{SRS}$) is to ignore the strength of the weak wave,
while preserving the strength of the strong vortex sheet/entropy wave,
and to place the error in the non-physical wave in the following manner:

\smallskip
{\it Case $1$. A weak wave $\lbrace U_-, U_1 \rbrace$ collides with the strong
vortex sheet/entropy wave $\lbrace U_1, U_+ \rbrace$ from below}.
The Riemann problem $\lbrace U_-, U_+ \rbrace$ is solved as follows:
$$
\begin{cases}
U_{-} &\quad \mbox{for $\frac{y}{x} < \chi(U_1, U_+)$},\\
U_{2} &\quad \mbox{for $\chi(U_1, U_+) < \frac{y}{x} < \hat{\lambda}$},\\
U_{+} &\quad \mbox{for $\frac{y}{x} > \hat{\lambda}$},
\end{cases}
$$
with $\chi(U_1, U_+)$ as the speed of the strong vortex sheet/entropy wave,
and state $U_2$ is solved in a way that $\lbrace U_-, U_2 \rbrace$
is the strong vortex sheet/entropy wave starting from $U_-$ and $\chi(U_1, U_+) = \chi(U_-, U_2)$.
Hence, ({\bf SRS}) keeps the same strength of the strong vortex sheet/entropy wave,
and the error appears in the non-physical fronts.

\smallskip
{\it Case $2$. A weak wave $\lbrace U_2, U_+ \rbrace$ collides with the strong vortex
sheet/entropy wave $\lbrace U_-, U_2 \rbrace$ from above}.
The Riemann problem $\lbrace U_-, U_+ \rbrace$ is solved as follows:
$$
\begin{cases}
U_{-} &\quad \mbox{for $\frac{y}{x} < \chi(U_-, U_2)$},\\
U_{2} &\quad \mbox{for $\chi(U_-, U_2) < \frac{y}{x} < \hat{\lambda}$},\\
U_{+} &\quad \mbox{for $\frac{y}{x} > \hat{\lambda}$},
\end{cases}
$$
with $\chi(U_-, U_2)$ denoting the speed of the strong vortex sheet/entropy wave.
\vspace{4mm}

\noindent\emph{3.2. Construction of Wave Front Tracking Approximations}
\vspace{2mm}

Given $\vartheta>0$, the corresponding front tracking approximate solution $U^{\vartheta}(x, y)$
is constructed as follows:
At $x = 0$, all the Riemann problems in $\overline{U}^{\vartheta}$ are solved by using
the accurate Riemann solver.
Furthermore, we can change the speed of one of the incoming fronts so that,
at any time $x > 0$, there is at most one collision involving only two incoming fronts.
This adjustment of speed can be chosen arbitrarily small.
Let $\omega_{\vartheta}$ be a fixed small parameter with $\omega_{\vartheta} \rightarrow 0$
as $\vartheta \rightarrow 0$, which will be determined later.
For convenience, subscript $j$ in $\alpha_{j}$ will be dropped henceforward,
and we will write $\alpha_{j}$ as $\alpha$ when no ambiguity arises
and employ the same notation $\alpha$
as a wave and its strength as before; the same applies for $\beta$.

\smallskip
\emph{Case $1$. Two weak waves with strengths $\alpha$ and $\beta$ interact
at some $x > 0$}.  The Riemann problem produced by this collision is solved in the following way:

\begin{itemize}
\item If $|\alpha \beta| > \omega_{\vartheta}$ and the two waves are physical,
then the accurate Riemann solver is employed.
\item If $|\alpha \beta| < \omega_{\vartheta}$ and the two waves are physical,
or there is a non-physical wave, then the simplified Riemann solver is employed.
\end{itemize}

\emph{Case $2$. A weak wave $\alpha$ interacts with the strong vortex sheet/entropy wave
and one weak wave at some $x > 0$}. The Riemann problem produced by this collision is solved in the following way:
\begin{itemize}
\item If $|\alpha| > \omega_{\vartheta}$ and the weak wave is physical, then  the accurate Riemann solver is applied.
\item If $|\alpha| < \omega_{\vartheta}$ and the weak wave is physical, or this wave is non-physical,
then the simplified Riemann solver is applied.
\end{itemize}

\emph{Case $3$. The flow perturbation due to the Lipschitz wall boundary}.
\begin{itemize}
\item When the change of angle of the boundary is larger than $\omega_{\vartheta}$
and the weak wave is physical, then the accurate Riemann solver is employed to solve the lateral Riemann problem.
\item If the change of angle of the boundary is less than $\omega_{\vartheta}$, then this perturbation is ignored.
\end{itemize}

\emph{Case $4$. The physical wave collides with the boundary}.
The accurate Riemann solver is employed to solve the lateral Riemann problem.

\vspace{3mm}
\noindent
\emph{3.3. Glimm's Functional and Wave Interaction Potential}
\vspace{2mm}

The goal in this subsection is to construct the suitable Glimm-type functional and
the associated wave interaction potential $\mathcal{Q}$ for the initial-boundary value problem
\eqref{2.1}--\eqref{1.8}.
This involves a careful combination
of the additional nonlinear waves generated
from the wall boundary vertices.
\vspace{3mm}

$\textbf{Definition 3.1}$ (\emph{Approaching waves}).

(i) Two weak fronts $\alpha$ and $\beta$,
located at points $y_{\alpha} < y_{\beta}$ and of the characteristic
families $j_{\alpha}$, $j_{\beta}$ $\in$ $\lbrace 1, \ldots, 4\rbrace$, respectively, are said to be
approaching each other if the following two conditions are concurrently satisfied:

\begin{itemize}
\item $y_{\alpha}$ and $y_{\beta}$ are both in one of the two intervals into which $\mathbb{R}$
is partitioned by the location of the strong vortex sheet/entropy wave.
That is, both waves are either in $\Omega_{-}$ or $\Omega_{+}$;
\item Either $j_{\alpha} > j_{\beta}$ or else $j_{\alpha}$ = $j_{\beta}$
and at least one of them is a shock.
\end{itemize}

\noindent In this case, we write $(\alpha, \beta)$ $\in$ $\mathcal{A}$.

\smallskip
(ii) A weak wave $\alpha$ of the characteristic family $j_{\alpha}$ is said to be approaching
the strong vortex sheet/entropy wave if either $\alpha \in \Omega_{-}$ and $j_{\alpha} = 4$,
or $\alpha \in \Omega_{+}$ and $j_{\alpha} = 1$. We then write $\alpha \in \mathcal{A}_{v/e}$.

\smallskip
(iii) A weak wave $\alpha$ of the characteristic family $j_{\alpha}$
is said to be approaching the boundary if $\alpha \in \Omega_{-}$ and $j_{\alpha} = 1$.
We then write $\alpha \in \mathcal{A}_{b}$.
\vspace{4mm}

Define the total (weighted) strength of weak waves in $U^{\vartheta}(x, \cdot)$ as
$$
\mathcal{V}(x) = \sum_{\alpha} |b_{\alpha}|,
$$
where, for a weak wave $\alpha$ of the $j$-family, its weighted strength is defined as
\begin{equation}\label{3.1}
b_{\alpha} =
\begin{cases} k_{+}\alpha &\quad \text{if $\alpha \in \Omega_{+}$ and $j_{\alpha} =1$},
\\
\alpha &\quad \text{if $\alpha \in \Omega_{-}$},
\end{cases}
\end{equation}
where $k_{+} = \frac{{2}K_{21}}{K^{*}}$, and coefficient $K_{21}$ is given as in Lemma 2.6.

Next, the wave interaction potential $\mathcal{Q}(x)$ is defined as
\begin{align}
\mathcal{Q}(x) &= C^{\ast} \sum_{(\alpha, \beta) \in \mathcal{A}} |b_{\alpha} b_{\beta}|
+ K^{\ast} \sum_{\alpha \in \mathcal{A}_{v/e}} |b_{\alpha}| + \sum_{\beta \in \mathcal{A}_{b}} |b_{\beta}|
+ \widetilde{K_{b0}} \sum_{a_{l} > x} |\theta_l|
\nonumber \\
&= \mathcal{Q}_{\mathcal{A}} + \mathcal{Q}_{v/e} + \mathcal{Q}_{b} + \mathcal{Q}_{\Theta}, \label{3.2}
\end{align}
where $K^{\ast} \in (K_{11}, 1)$ and $\widetilde{K_{b0}} > K_{b0}$,
while $C^{\ast}$ is a constant to be specified later.
To control the total variation of the new waves produced by the boundary vertices,
$\mathcal{Q}_{\Theta}$ in our wave interaction potential $\mathcal{Q}(x)$ is an added term,
compared to that for the Cauchy problem.

\smallskip
\emph{The Glimm-type functional $\mathcal{G}$ is defined as follows}:
\begin{equation}
\mathcal{G}(x)
= \mathcal{V}(x) + \Ka \mathcal{Q}(x)
  + |U^{\diamond}(x) - U^{+}_{0}| + |U_{\diamond}(x) - U^{-}_{0}|,
\end{equation}
\emph{where
states $U_{\diamond}(x)$ and $U^{\diamond}(x)$ are
the below state and the above state of
the strong vortex sheet/entropy wave at {\it time} $x$
respectively,
$U^{-}_{0}$ and $U^{+}_{0}$
are the below and above state
of the strong vortex sheet/entropy wave at $x = 0$
respectively,  and $\Ka$ is a large positive constant to be determined later.}

\smallskip
Notice that $\mathcal{V}$, $\mathcal{Q}$, and $\mathcal{G}$ remain unchanged between
any pair of subsequent interaction times.
However, we will show that, across an interaction {\it time} $x$,
both $\mathcal{Q}$ and $\mathcal{G}$ decrease.

\medskip
\noindent $\textbf{Lemma 3.1.}$
\emph{Assume that ${\rm TV}(\widetilde{U_0}(\cdot)) + {\rm TV}(g^{\prime}\left(\cdot\right))$ is
sufficiently small.
Then $\mathcal{V}(x)$ remains sufficiently small for all $x > 0$,
and ${\rm TV}(U^{\vartheta}(x, \cdot))$ has a uniform bound
for any $\vartheta > 0$}.
\vspace{3mm}

\noindent
$\textbf{Proof.}$
With the Glimm-type functional $\mathcal{G}$, consider
$$
\Delta \mathcal{G}(x) = \mathcal{G}(x^{+}) - \mathcal{G}(x^{-}),
$$
where $x^{-}$ and $x^{+}$ denote the {\it times} before and after
the interaction {\it time} $x > 0$, respectively.

\smallskip
\emph{Case $1$. Two weak waves $\alpha$ and $\beta$ collide}.
States $U^{\diamond}\left(x\right)$ and $U_{\diamond}\left(x\right)$
do not alter across this interaction {\it time} $x > 0$, so that
\begin{align}
\Delta \mathcal{G}(x)
&= \mathcal{V}(x^{+}) - \mathcal{V}(x^{-}) + \Ka\left(\mathcal{Q}(x^{+}) - \mathcal{Q}(x^{-})\right)     \nonumber \\
&\le \mathcal{B}_{1}|b_{\alpha}b_{\beta}| - \Ka \left((C^{\ast}-\mathcal{B}_{0})|b_{\alpha}b_{\beta}|
    - C^{\ast} |b_{\alpha}b_{\beta}|\mathcal{V}(x^{-})\right), \nonumber
\end{align}
where $\mathcal{B}_{0}$ and $\mathcal{B}_{1}$ are constants independent of $\vartheta$.

\smallskip
\emph{Case $2$. A weak wave $\alpha$ of the $1$-family interacts with the boundary}.
$$
\Delta \mathcal{G}(x)
= K_{b1}\alpha - \alpha - \Ka\left((1-K^{\ast}K_{b1})\alpha- C^{*}K_{b1}\mathcal{V}( x^{-})\alpha \right),
$$
where $K^{\ast}K_{b1} < 1$.

\smallskip
\emph{Case $3$. A new $4$-wave $\alpha$ produced by the Lipschitz wall boundary}.
$$
\Delta \mathcal{G}(x)
= K_{b0}\theta_{l}- \Ka\left((\widetilde{K_{b0}}-K^{*}K_{b0})\theta_{l}
  -C^{\ast}K_{b0}\theta_{l}\mathcal{V}( x^{-})\right),
$$
where $K_{b0} < \widetilde{K_{b0}}$ is large.

\medskip
In the following two cases,
states $U_{\diamond}(x)$ and $U^{\diamond}(x)$
change across this interaction {\it time} $x > 0$.

\smallskip
\emph{Case $4$. A weak wave $\alpha$ of the $4$-family collides with the strong vortex sheet/entropy wave from below}.
\begin{align*}
 \Delta \mathcal{G}(x)
 &\le \mathcal{V}(x^{+}) - \mathcal{V}(x^{-})
   + |U^{\diamond}( x^{+}) - U^{\diamond}( x^{-})|
   + |U_{\diamond}(x^{+}) - U_{\diamond}( x^{-})|
   + \Ka\left(\mathcal{Q}(x^{+}) - \mathcal{Q}(x^{-})\right)  \\
   &\le  (K_{11}+K_{14}+C^*)\alpha - \alpha
       - \Ka\Big((K^{\ast}-K_{11})\alpha - C^{\ast}(K_{11}+K_{14})\mathcal{V}(x^{-})\alpha\Big).
\end{align*}

\emph{Case $5$. A weak wave $\alpha$ of the $1$-family collides with the strong vortex sheet/entropy wave from above}.
$$
\Delta \mathcal{G}(x)
\le (K_{21}+K_{24}+C^*)\alpha - b_{\alpha}
  - \Ka\Big( (K^{\ast}b_{\alpha}- K_{21}\alpha)
      -C^{\ast}( K_{21}+K_{24})\mathcal{V}(x^{-})\alpha \Big).
$$

In the cases above, $K_{11} < K^{\ast} < 1$, $K^* b_{\alpha} \ge {2}K_{21}|\alpha|$ in connection with weight $k_{+}$,
and
the constant $C^{\ast} > \mathcal{B}_{0} > 0$ is large.

\medskip
Next, we establish that the {\it total $($weighted$)$ strength of waves}
in $U^{\vartheta}(x, \cdot)$ remains sufficiently small for all $x > 0$
if it is sufficiently small at $x=0$. More precisely,
$$
\mathcal{V}(x)\ll 1 \qquad\,\, \text{  for all $x > 0$}.
$$
This can be proved as follows:

\vspace{3mm}
\emph{{\rm (i)}  $x_1 > 0$ is the first interaction time}.
Given that $\mathcal{V}( x^{-}_{1})=\mathcal{V}(0)\leq {\rm TV}(\widetilde{U_0}(\cdot)) \ll 1$
and $\sum_{\l = 0}^{\infty}\theta_{l} \leq {\rm TV}(g^{\prime}(\cdot)) \ll 1$ in Cases 1--5 above,
we conclude that, for $\Ka$ sufficiently large and $\omega_{\vartheta}$ sufficiently small,
$$
\Delta \mathcal{G}(x_{1}) \leq 0, \qquad {i.e.,}
\quad \mathcal{G}(x_{1}^{+}) \leq \mathcal{G}(x_{1}^{-}) = \mathcal{G}(0).
$$
Therefore,
\begin{align}
\mathcal{V}(x_{1}^{+})
&\leq \mathcal{G}(x_{1}^{+}) \leq \mathcal{G}(0) \leq \mathcal{V}(0) + \Ka \mathcal{Q}(0) \nonumber \\
&= \mathcal{V}(0) + \Ka\Big(C^{\ast}\mathcal{V}^{2}(0) + \mathcal{V}(0)
+ \widetilde{K_{b0}}\sum_{l = 0}^{\infty}\theta_{l}\Big)  \nonumber \\
&\le C\Big(\mathcal{V}(0) + \sum_{l = 0}^{\infty}\theta_{l}\Big) \ll 1. \nonumber
\end{align}

\emph{{\rm (ii)} $\mathcal{V}(x_{m}^{-}) \ll 1$ and $\mathcal{G}(x_{m}^{+}) \leq \mathcal{G}(x_{m}^{-})$ for any $m < n$}.
Then, for the next interaction {\it time} $x_{n}$, similar to Case (i), we also conclude
$$
\Delta \mathcal{G}(x_{n}) \leq 0, \qquad {i.e.,} \quad  \mathcal{G}(x_{n}^{+}) \leq \mathcal{G}(x_{n}^{-}) = \mathcal{G}(x_{n-1}^{+}).
$$

Therefore, all together, we obtain
\begin{eqnarray}
&&\mathcal{V}(x_{n}^{+})
 +|U^{\diamond}(x_{n}^{+}) - U_{0}^{+}| + |U_{\diamond}(x_{n}^{+}) - U_{0}^{-}|
\nonumber\\[1mm]
&& \leq \mathcal{G}(x_{n}^{+}) \leq \mathcal{G}(x_{n}^{-}) = \mathcal{G}(x_{n-1}^{+}) \leq \ldots \leq \mathcal{G}(0)\nonumber \\[1mm]
&&= \mathcal{V}(0) + \Ka \mathcal{Q}(0)     \nonumber \\
&&= \mathcal{V}(0) + \Ka\Big(C^{\ast}\mathcal{V}^{2}(0) + \mathcal{V}(0) + \widetilde{K_{b0}}\sum_{l = 0}^{\infty}\theta_{l}\Big)  \nonumber \\
&&\le C\Big(\mathcal{V}(0) + \sum_{l = 0}^{\infty}\theta_{l}\Big) \ll 1. \nonumber
\end{eqnarray}
This implies that $\mathcal{V}(x) \ll 1$ for all $x > 0$, since $C$ is independent of $x$.

Furthermore, the total variation of $U^{\vartheta}(x, \cdot)$ is uniformly bounded:
\begin{equation}\label{3.4}
{\rm TV}\lbrace U^{\vartheta}(x, \cdot)\rbrace
\approx \mathcal{V}(x)
|U^{\diamond}(x) - U_{0}^{+}| + |U_{\diamond}(x) - U_{0}^{-}|
+ |\sigma_{20}| + |\sigma_{30}| = \mathcal{O}(1).
\end{equation}
This completes the proof.

\medskip
In order to define the front tracking approximate solution $U^{\vartheta}(x,\cdot)$ for any $x > 0$,
along with a uniform bound on the total variation,
we also need to have a finite number of wave-fronts in $U^{\vartheta}(x,\cdot)$.
This is given by the following lemma.

\medskip
\noindent $\textbf{Lemma 3.2.}$ \emph{For any fixed $\vartheta > 0$ small enough,
the number of wave-fronts in $U^{\vartheta}(x, y)$ is finite
and the approximate solutions $U^{\vartheta}(x, y)$ are defined for all $x > 0$.
Moreover, for any $x > 0$, the total strength of all the non-physical waves is
 of order $\mathcal{O}(1)\left(\delta_{\vartheta} + \omega_{\vartheta}\right)$}.

\medskip
\noindent $\textbf{Proof.}$ We first note that the total interaction potential $\mathcal{Q}(x)$ remains
unchanged when there is no interaction and
decreases across an interaction {\it time} $x > 0$, as discussed in Cases 1--5 in Lemma 3.1.
Furthermore, from Cases 1--5 and the subsequent analysis above,
we have concluded that $\mathcal{V}(x) \ll 1$. Thus, we can fix some number $\nu\in (0, 1)$ such that
\begin{align}\label{3.5}
 \Delta \mathcal{Q}(x)
 &= \mathcal{Q}(x^{+}) - \mathcal{Q}(x^{-})      \nonumber \\
 &\le
 \begin{cases}
 -\nu |b_{\alpha}b_{\beta}| &\mbox{ if both waves $\alpha$ and $\beta$ are weak,} \\
 -\nu |b_{\alpha}| &\mbox{ if the weak wave $\alpha$ hits the strong vortex sheet/entropy wave,} \qquad \\
 -\nu |\theta_{l}| &\mbox{ if the angle  of the boundary changes.}
        \end{cases}
\end{align}
Now, following an argument similar to the one given in \ci{Baiti-Jenssen-1998},
we reach the following conclusions:
Note that initially $\mathcal{Q}(0)$ is bounded and $\mathcal{Q}$ decreases thereafter for each case.
Moreover, in the case where the interaction potential between the incoming waves
or the change of angle of the boundary is larger than $\omega_{\vartheta}$,
$\mathcal{Q}$ decreases by at least $\nu \omega_{\vartheta}$ in these interactions,
as implied by the bounds given in \eqref{3.5}.
Following the wave-front tracking method in our problem,
new physical waves can only be produced by such interactions.
Furthermore, when the weak wave $\alpha$ of $1$-family collides with the wall boundary,
we have solved the lateral Riemann problem and shown that,
after this interaction, there is only a reflected wave of $4$-family with
the reflection coefficient $1$.
Thus, before and after this interaction,
the number of the waves stays the same, which implies that the number of the waves is finite.
Finally, because the non-physical waves are generated only when the physical waves collide,
we can also conclude that the number of non-physical wave fronts is finite;
if two waves can only collide once, the number of interactions is also finite.
Consequently, it follows that the approximate solutions $U^{\vartheta}(x, \cdot)$ are defined
for all $x > 0$.
The similar argument allows us to conclude that the total strength of all the non-physical wave fronts
at any $x$ is of order $\mathcal{O}(1)(\delta_{\vartheta}+\omega_{\vartheta})$.
This completes the proof.
\vspace{4mm}

Following the line of arguments as in \ci{Baiti-Jenssen-1998,Bressan-1992}
for the wave-front tracking algorithm and Lemma 3.1 above,
we conclude this section with the following theorem for the global existence
of entropy solutions of the initial-boundary value
problem \eqref{2.1}--\eqref{1.8}.
\vspace{4mm}

\noindent $\textbf{Theorem 3.1.}$ \emph{Suppose that ${\rm TV}(\widetilde{U_0}(\cdot))+{\rm TV}(g^{\prime}(\cdot))$
is suitably small. Then, for the initial-boundary value problem \eqref{2.1}--\eqref{1.8},
there exists a global entropy solution in {\rm BV} satisfying the {\it steady Clausius entropy}
inequality \eqref{1.10}.}

\section{\small{The Lyapunov Functional for the ${L}^1$--Distance between Two Solutions}}

To show that the wave-front tracking approximations, constructed for the existence analysis in \S 3,
converge to a unique limit, we estimate the distance between
any two $\vartheta$-approximate $U$ and $V$ of problem \eqref{2.1}--\eqref{1.8}.
To this end, we develop the Lyapunov functional $\Phi(U, V)$, equivalent to the $L^{1}$--distance:
\[
{C}^{-1} \, \norm{U(x,\cdot) - V(x,\cdot)}_{L^{1}}
\leq \mathit{\Phi}(U, V) \leq {C} \, \|U(x,\cdot) - V(x,\cdot)\|_{L^{1}},
\]
and prove that $\Phi(U,V)$ is {\it almost decreasing}:
\[
\mathit{\Phi} \left(U(x_{2},\cdot), V(x_{2},\cdot)\right)
- \mathit{\Phi}\left(U(x_{1},\cdot), V(x_{1},\cdot)\right)
\leq C\vartheta(x_{2} - x_{1}) \hspace{5mm} \text{for all } x_{2} > x_1 > 0,
\]
for some constant ${C}>0$.
Here $U$ and $V$ are two approximate solutions constructed via the wave-front tracking method,
and the small approximation parameter $\vartheta$ is responsible
for controlling the subsequent errors:

\begin{my_itemize}
\item Errors in the approximation of the initial data and the boundary.
\item Errors in the speeds of shocks, vortex sheets, entropy waves, and rarefaction fronts.
\item The total strength of non-physical fronts.
\item The maximum strength of rarefaction fronts.
\end{my_itemize}

\smallskip
Along the line of arguments presented in \ci{Bressan-Liu-Yang-1999, Lewicka-Trivisa-2002, Liu-Yang-1999},
with {\it time} $x$ fixed, at each $y$, one connects state $U(y)$ with $V(y)$ in the state space
by going along the Hugoniot curves $S_{1}, C_{2}, C_{3}$, and $S_{4}$.
Depending on the location of the strong vortex sheet/entropy wave in $U(y)$ and $V(y)$,
the distance between $U(y)$ and $V(y)$ is estimated along discontinuity waves in possibly different {\it directions},
determining the strength of the $j$-Hugoniot wave, $h_{j}(y)$, in the following way:
\begin{my_itemize}
\item If $U(y)$ and $V(y)$ are both in $\Omega_{-}$ and $\Omega_{+}$,
then it begins at state $U(y)$ and moves along the Hugoniot curves to reach state $V(y)$.
\item If $U(y)$ is in $\Omega_{-}$ and $V(y)$ is in $\Omega_{+}$, then it begins at state $U(y)$
  and moves along the Hugoniot curves to reach state $V(y)$.
\item If $V(y)$ is in $\Omega_{-}$ and $U(y)$ is in $\Omega_{+}$, then it begins at state $V(y)$
and moves along the Hugoniot curves to reach state $U(y)$.
\end{my_itemize}

\smallskip
\noindent Define the $L^{1}$--weighted strengths of the waves in the solution
of the Riemann problem $\{U(y), V(y)\}$ or $\{V(y), U(y)\}$ as follows:
\begin{equation}\label{4.1}
 q_{j}(y)=
 \begin{cases}
 w^{b}_{j}\,h_{j}(y)&\mbox{ whenever $U(y)$ and $V(y)$ are both in $\Omega_{-}$,} \\
 w^{m}_{j}\, h_{j}(y)&\mbox{ whenever $U(y)$ and $V(y)$ are both in different domains,} \\
 w^{a}_{j}\,h_{j}(y)&\mbox{ whenever $U(y)$ and $V(y)$ are both in $\Omega_{+}$,}
        \end{cases}
\end{equation}
with constants $w^{b}_{j}$, $w^{m}_{j}$, and $w^{a}_{j}$ above to be specified later on,
based on the estimates of wave interactions and reflections in Lemmas 2.2--2.7.

We define the following Lyapunov functional:
\begin{equation}\label{4.2}
\mathit{\Phi}(U, V) = \sum_{j = 1}^{4}\int\limits_{g(x)}^{\infty}\ |q_{j}(y)|W_{j}(y)\,{\rm d}y,
\end{equation}
\vspace{3mm}
where the weights are given by
\begin{equation}\label{4.3}
W_{j}(y) = 1 + \Ka_{1}A_{j}(y) + \Ka_{2}\big(\mathcal{Q}(U) + \mathcal{Q}(V)\big)
\end{equation}
with constants $\Ka_{1}$ and $\Ka_{2}$ to be determined later. Here $\mathcal{Q}$
denotes the total wave interaction potential incorporating
the boundary effect as defined in \eqref{3.2},
and $A_{j}(y)$ denotes the total strength of waves in $U$ and $V$, which approach
the $j$-wave $q_{j}(y)$, defined in the following
manner (for $y$ where there is no jump in $U$ or $V$):
\begin{equation}\label{4.4}
A_{j}(y) = F_{j}(y) + G_{j}(y)
+ \begin{cases}
H_{j}(y)&\mbox{ if $j$-wave $q_{j}(y)$ is small and the $j$-field is genuinely nonlinear,}\\[1mm]
0&\mbox{ if $j = 2, 3$, and $q_{j}(y) = B$ is large.}
        \end{cases}
\end{equation}
We first define the following global weights $G_{j}$:

\bigskip
\begin{tabular}{ l | l || l || l }
\hline	
$G_{j}(y)\,$ & $U, V$ are both in $\Omega_{-}$ & $U, V$ are in distinct regions & $U, V$ are both in $\Omega_{+}$ \\		
\hline
$G_{1}(y)$ & 4B & 2B & 4B \\
$G_{2,3}(y)$ & 0 & 0 & 0 \\
$G_{4}(y)$ & 4B & 2B & 2B \\
\hline
\end{tabular}

\bigskip
\noindent
The summands in \eqref{4.4} are defined as follows:
\begin{align}
F_{j}(y) &= \sum_{\substack{\alpha \in \mathcal{J} \setminus \mathcal{SC}
  \\ y_{\alpha} < y, j < k_{\alpha} \leq 4}}|\alpha| + \sum_{\substack{\alpha \in \mathcal{J} \setminus \mathcal{SC}
  \\ y_{\alpha} > y, 1 \leq k_{\alpha} < j}}
   |\alpha |,       \nonumber \\[2mm]
H_{j}(y) &= \begin{cases}
(\sum_{\substack{\alpha \in \mathcal{J}(U) \setminus \mathcal{SC}, y_{\alpha} < y, k_{\alpha} = j}}
+ \sum_{\substack{\alpha \in \mathcal{J}(V) \setminus \mathcal{SC}, y_{\alpha} > y, k_{\alpha} = j}})|\alpha|
&\quad \mbox{ if $q_{j}(y) < 0$,} \\[1mm]
(\sum_{\substack{\alpha \in \mathcal{J}(V) \setminus \mathcal{SC}, y_{\alpha} < y, k_{\alpha} = j}}
+ \sum_{\substack{\alpha \in \mathcal{J}(U) \setminus \mathcal{SC}, y_{\alpha} > y, k_{\alpha} = j}})|\alpha|
&\quad \mbox{ if $q_{j}(y) > 0$,}
\end{cases}       \nonumber
\end{align}
where, at each $x$, $\alpha$ stands for the (non-weighted) strength of wave $\alpha \in \mathcal{J}$,
located at point $y_{\alpha}$ and belonging to the characteristic family $k_{\alpha}$;
$\mathcal{J} = \mathcal{J}(U) \cup \mathcal{J}(V)$ and $\mathcal{SC} = \mathcal{SC}(U) \cup \mathcal{SC}(V)$
are the set of all the waves (in $U$ and $V$) and the set of all the strong characteristic
discontinuities (in $U$ and $V$), respectively.

\medskip
\noindent Under the assumption that
${\rm TV}(\widetilde{U_0}(\cdot)) + {\rm TV}(\widetilde{V_0}(\cdot)) + {\rm TV}(g^{\prime}(\cdot))$
is small enough with $U(x, \cdot)$, $V(x, \cdot)$ $\in {\rm BV} \cap L^{1}$, one concludes
\begin{align}
&\mathcal{M}^{-1} \norm{U(x,\cdot) - V(x,\cdot)}_{L^{1}}
\leq \sum_{j=1}^{4} \int\limits^{\infty}\limits_{g(x)} |q_{j}(y)|\,{\rm d}y
\leq \mathcal{M} \|U(x,\cdot) - V(x,\cdot)\|_{L^{1}}, \nonumber \\
&1 \leq W_{j}(y) \leq \mathcal{M}, \hspace{2mm} j = 1, \ldots, 4, \nonumber
\end{align}
where constant $\mathcal{M}$ is independent of $\vartheta$ and {\it time} $x$.
Here \emph{we define the strength of any large wave of the $2$-characteristic or $3$-characteristic family
to equal to some fixed number B $($larger than all the strengths of the small waves$)$,}
and the terms ``small" and ``large" refer to the waves
that connect the states in the same or in the distinct domains $\Omega^{-}$ and $\Omega^{+}$,
respectively.

\smallskip
Consequently, we have
\begin{equation}\label{4.5}
{C}^{-1} \norm{U(x,\cdot) - V(x,\cdot)}_{L^{1}}
\leq \mathit{\Phi}(U, V) \leq {C} \|U(x,\cdot) - V(x,\cdot)\|_{L^{1}}
\end{equation}
for any $x \geq 0$ with constant ${C}>0$ depending only on the quantities
independent of $x$: the strength of the strong vortex sheet/entropy wave and
${\rm TV}(\widetilde{U_0}(\cdot)) + {\rm TV}(\widetilde{V_0}(\cdot)) + {\rm TV}(g^{\prime}(\cdot))$.

\section{The $L^{1}$--Stability Estimates}
In this section, we establish the $L^1$--stability estimates.

\vspace{3mm}
\noindent
\emph{5.1. Evolution of the Lyapunov Functional $\mathit{\Phi}$ in the Flow Direction $x > 0$}.

\vspace{2mm}
For each $j=1, \ldots, 4$, $\lambda_{j}(y)$ is the speed of the $j$-wave $q_{j}(y)$
(along the Hugoniot curve in the phase space).
Then, at a {\it time} $x > 0$ that is not the interaction time of the waves
in either $U(x)=U(x,\cdot)$ or $V(x)=V(x,\cdot)$,
an explicit computation gives
\begin{align}
& {{\rm d}\over {\rm d}x}\mathit{\Phi}(U(x), V(x)) \nonumber \\
&= \sum_{\alpha \in \mathcal{J}}\sum_{j = 1}^{4}\left(\abs{q_{j}(y_{\alpha}^{-})}W_{j}(y_{\alpha}^{-})
      - \abs{q_{j}(y_{\alpha}^{+})}W_{j}(y_{\alpha}^{+})\right)\dot{y}_{\alpha}
      - \sum_{j = 1}^{4}\abs{q_{j}(b)}W_{j}(b)\dot{y}_{b}    \nonumber \\
&= \sum_{\alpha \in \mathcal{J}}\sum_{j = 1}^{4}\left(\abs{q_{j}(y_{\alpha}^{-})}W_{j}(y_{\alpha}^{-})\big(\dot{y}_{\alpha} - \lambda_{j}(y_{\alpha}^{-})\big) - \abs{q_{j}(y_{\alpha}^{+})}W_{j}(y_{\alpha}^{+})\big(\dot{y}_{\alpha} - \lambda_{j}(y_{\alpha}^{+})\big)\right)     \nonumber \\
& \quad + \sum_{j = 1}^{4}\abs{q_{j}(b)}W_{j}(b)\big(-\dot{y}_{b} + \lambda_{j}(b)\big), \label{4.6}
\end{align}
where $\dot{y}_{\alpha}$ denotes the speed of the Hugoniot wave $\alpha \in \mathcal{J}$,
$b = g(x)^{+}$ stands for the points close to the boundary, and $\dot{y}_{b}$ is the slope of the boundary.

Then \eqref{4.6} can be written as
\begin{equation}\label{4.9}
{{\rm d}\over {\rm d}x}\mathit{\Phi}\left(U(x), V(x)\right)
= \sum_{\alpha \in \mathcal{J}}\sum_{j=1}^{4}E_{\alpha,j} + \sum_{j=1}^{4}E_{b,j},
\end{equation}
where
\begin{align}
E_{\alpha,j}& =\abs{q_{j}^{+}}W_{j}^{+}\big(\lambda_{j}^{+} - \dot{y}_{\alpha} \big)
- \abs{q_{j}^{-}}W_{j}^{-}\big(\lambda_{j}^{-} - \dot{y}_{\alpha} \big), \label{4.7}\\[2mm]
E_{b,j}& = \abs{q_{j}(b)}W_{j}(b)\big(-\dot{y}_{b} + \lambda_{j}(b)\big), \label{4.8}
\end{align}
with $q_{j}^{\pm}$ = $q_{j}(y^{\pm}_{\alpha})$, $W^{\pm}_{j} = W_{j}(y^{\pm}_{\alpha})$,
and $\lambda^{\pm}_{j}$ = $\lambda_{j}(y^{\pm}_{\alpha})$.

\smallskip
Our central aim in \S 5.2 below is to prove the bounds:
\begin{align}
& \sum_{j=1}^{4}E_{\alpha,j} \leq \mathcal{O}(1)\vartheta\abs{\alpha}\quad
 \text{ when $\alpha$ is a weak wave in $\mathcal{J}$}, \label{4.10} \\
& \sum_{j=1}^{4}E_{\alpha,j} \leq \mathcal{O}(1)\abs{\alpha}
\quad  \text{ when $\alpha$ is a non-physical wave in $\mathcal{J}$,}\label{4.11}\\
& \sum_{j=1}^{4}E_{\alpha,j}
  \leq \mathcal{O}(1)B\vartheta
\quad \text{ when $\alpha$ is a strong vortex sheet/entropy wave in $\mathcal{J}$,} \label{4.12} \\
& \sum_{j=1}^{4}E_{b,j} \leq 0 \quad \text{ near the boundary}, \label{4.13}
\end{align}
where the quantities denoted by the Landau symbol $\mathcal{O}$(1) are independent of
constants $\Ka_1$ and $\Ka_2$.

With these bounds \eqref{4.10}--\eqref{4.13} together with the uniform bound on the total strengths
of waves \eqref{3.4}, we obtain
\begin{equation}\label{4.14}
{{\rm d}\over {\rm d}x}\mathit{\Phi}(U(x), V(x)) \leq \mathcal{O}(1)\vartheta.
\end{equation}
Integration of \eqref{4.14} over the interval $\left[0, x\right]$ yields
\begin{equation}\label{4.14a}
\mathit{\Phi}\left(U(x), V(x)\right) \leq \mathit{\Phi}\left(U(0), V(0)\right) + \mathcal{O}(1)\vartheta x.
\end{equation}

We remark that, at each interaction {\it time} $x$ when two fronts of $U$ or two fronts of $V$ interact,
by the Glimm interaction estimates, all the weight functions $W_{j}(y)$ decrease,
if constant $\Ka_{2}$ in the Lyapunov functional is taken to be sufficiently large.
Furthermore, due to the self-similar property of the Riemann solutions,
$\mathit{\Phi}$ decreases at this {\it time}.

\vspace{3mm}
\noindent
\emph{5.2. {Estimates for Bounds \eqref{4.10}--\eqref{4.13}}}.

\vspace{2mm}
We now establish
bounds \eqref{4.10}--\eqref{4.13}, particularly \eqref{4.12}--\eqref{4.13},
when $\alpha$ is a strong vortex sheet/entropy wave in $\mathcal{J}$ and
near the Lipschitz wall boundary, respectively.

\vspace{2mm}
For the case that the weak wave
$\alpha \in \mathcal{J} \mathrel{\mathop:}= \mathcal{J}(U) \cup \mathcal{J}(V)$
and the non-physical waves in $\mathcal{J}$,
which appears when $U$ and $V$ are both in $\Omega_{-}$ or $\Omega_{+}$,
estimates \eqref{4.10}--\eqref{4.11} are shown similarly based on
the arguments in Bressan-Liu-Yang \ci{Bressan-Liu-Yang-1999},
provided that $\frac{2|B|}{|\sigma_{20}|+|\sigma_{30}|}$ is
sufficiently small and $\Ka_{1}$ is sufficiently large.
In what follows, we focus only on the other two cases, namely \eqref{4.12}--\eqref{4.13}.

\smallskip
{\it Case $1$}. {\it The first strong vortex sheet/entropy wave $\alpha$ in $U$ or $V$ is crossed}.
Using Lemma 2.4, we have the estimates:
\begin{align}
h^{+}_{1} &= h^{-}_{1} + K_{11}h^{-}_{4}, \label{5.1} \\
h^{+}_{4} &= K_{14}h^{-}_{4}. \label{5.2}
\end{align}

\noindent Moreover, the essential estimate $|K_{11}| < 1$ given in Lemma 2.4 ensures
the existence of desired weights $w^{b}_{1}$ and $w^{b}_{4}$ in the following way.
\vspace{3mm}

\noindent
$\textbf{Lemma 5.1.}$
\emph{There exist $w^{b}_{1}$, $w^{b}_{4}$,
and $\gamma_{b}$ such that}
\begin{eqnarray}
&&\frac{w^{b}_{4}}{w^{b}_{1}} < 1, \label{5.3} \\
&&\frac{w^{b}_{1}}{w^{b}_{4}}\left|\frac{\lambda_{1}^{-} - \lambda_{2,3}}{\lambda_{4}^{-} - \lambda_{2,3}}\right|K_{11}
< \gamma_{b} < 1. \label{5.4}
\end{eqnarray}
\qquad

With Lemma 5.1, we estimate $E_{j}$ for $j = 1, \ldots, 4$, starting with $E_1$:
By \eqref{5.1} and \eqref{5.4},
\begin{align*}
E_{1} &= |q_{1}^{-}|(\lambda_{1}^{-}-\dot{y}_{\alpha})(W_{1}^{+} - W_{1}^{-})
         + W^{+}_{1}\left(|q^{+}_{1}|(\lambda_{1}^{+}-\dot{y}_{\alpha})
      -|q^{-}_{1}|(\lambda^{-}_{1}-\dot{y}_{\alpha})\right)          \nonumber \\
&= 2B\Ka_{1}w^{b}_{1}|h^{-}_{1}||\lambda_{1}^{-}-\dot{y}_{\alpha}|
   + W^{+}_{1}\left(|q^{-}_{1}||\lambda_{1}^{-}-\dot{y}_{\alpha}|-|q^{+}_{1}||\lambda_{1}^{+}
    -\dot{y}_{\alpha}|\right) \nonumber \\
&\le 2B\Ka_{1}w^{b}_{1}|h^{-}_{1}||\lambda_{1}^{-}-\dot{y}_{\alpha}|
  + W^{+}_{1}\big(w^{b}_{1}|h^{+}_{1}||\lambda_{1}^{-}-\dot{y}_{\alpha}|
   + w^{b}_{1}K_{11}|h^{-}_{4}||\lambda_{1}^{-}-\dot{y}_{\alpha}|
    - w^{m}_{1}|h^{+}_{1}||\lambda_{1}^{+}-\dot{y}_{\alpha}|\big) \nonumber \\
&\le 2B\Ka_{1}w^{b}_{1}|h^{-}_{1}||\lambda_{1}^{-}-\dot{y}_{\alpha}|
   + 2B\Ka_{1}\big(w^{b}_{1}|h^{+}_{1}||\lambda_{1}^{-}-\dot{y}_{\alpha}|
      + \gamma_{b}w^{b}_{4}|h^{-}_{4}|(\lambda_{4}^{-}-\dot{y}_{\alpha})
- w^{m}_{1}|h^{+}_{1}||\lambda_{1}^{+}-\dot{y}_{\alpha}|\big) \nonumber \\
&\quad + (\Ka_{1}A_{W^{+}_{1}} + \widehat{\Ka})|q^{-}_{1}||\lambda_{1}^{-}-\dot{y}_{\alpha}|
   - (\Ka_{1}A_{W^{+}_{1}} + \widehat{\Ka})|q^{+}_{1}||\lambda_{1}^{+}-\dot{y}_{\alpha}|, \nonumber
\end{align*}
where
\begin{equation}\label{5.5a}
\widehat{\Ka} := 1 + \Ka_{2}(\mathcal{Q}(U) + \mathcal{Q}(V))>0,
\end{equation}
$W^{+}_{1} = W_{1}(y_{\alpha}^+) = 2B\Ka_{1} + \Ka_{1}A_{W^{+}_{1}} + \widehat{\Ka}$,
$A_{W^{+}_{1}} = F_{1}(y_{\alpha}^+) + H_{1}(y_{\alpha}^+)$ here is the total strength
of all the weak waves in $U$ and $V$ which approach the $1$-wave $q_{1}^{+}=q_{1}(y_{\alpha}^+)$,
and $2B\Ka_{1}$ is from weight $G_{1}(y_{\alpha}^+)$.

\smallskip
For $j = 2, 3$, $W^{+}_{j} = W^{-}_{j}$ so that \eqref{4.6} reduces to
\begin{eqnarray*}
E_{j}= W^{-}_{j}\big(|q_{j}^{+}|(\lambda_{j}^{+}-\dot{y}_{\alpha}) - |q_{j}^{-}|(\lambda_{j}^{-}-\dot{y}_{\alpha})\big)
\le \mathcal{O}(1)B\Big(\vartheta + \sum_{i=1,4}|q_{i}^{-}| \Big),
\end{eqnarray*}
where $k \notin \lbrace j, 1, 4 \rbrace$.

\smallskip
For $j = 4$,
\begin{align}
E_{4} &= |q_{4}^{-}|(\lambda_{4}^{-}-\dot{y}_{\alpha})(W_{4}^{+} - W_{4}^{-}) + W^{+}_{4}\left(|q^{+}_{4}|(\lambda_{4}^{+}-\dot{y}_{\alpha})
-|q^{-}_{4}|(\lambda^{-}_{4}-\dot{y}_{\alpha})\right)          \nonumber \\
&= -2B\Ka_{1}|q^{-}_{4}|(\lambda_{4}^{-}-\dot{y}_{\alpha})
+ \big(2B\Ka_{1} + \Ka_{1}A_{W^{+}_{4}} + \widehat{\Ka}\big)|q^{+}_{4}|(\lambda_{4}^{+}-\dot{y}_{\alpha}) \nonumber \\
&\quad -{}\big(2B\Ka_{1} + \Ka_{1}A_{W^{+}_{4}} + \widehat{\Ka}\big)|q^{-}_{4}|(\lambda_{4}^{-}-\dot{y}_{\alpha}),  \nonumber
\end{align}
where $W^{+}_{4} = W_{4}(y_{\alpha}^+) = 2B\Ka_{1} + \Ka_{1}A_{W^{+}_{4}} + \widehat{\Ka}$ with constant $\widehat{\Ka}$ determined by \eqref{5.5a},
$A_{W^{+}_{4}} = F_{4}(y_{\alpha}^+) + H_{4}(y_{\alpha}^+)$ is the total strength
of all the weak waves in $U$ and $V$ which approach the $4$-wave $q_{4}^{+}=q_{4}(y_{\alpha}^+)$,
and $2B\Ka_{1}$ is from weight $G_{4}(y_{\alpha}^+)$.

For the weighted $L^{1}$--strength $q_{j}(y)$ in \eqref{4.1},
we choose $w^{b}_{1}$ to be small enough relative to $w^{m}_{1}$, $w^{b}_{4}$ large enough
relative to $w^{m}_{4}$,
and $\Ka_{1}$ large enough and the total variation of $U$ and $V$ small enough.
Then we use \eqref{5.1}--\eqref{5.2} to obtain
\begin{align*}
\sum^{4}_{j=1}E_{j}
\leq &\, 2B\Ka_{1}\big(w^{b}_{1}|h^{+}_{1}||\lambda_{1}^{-}-\dot{y}_{\alpha}|
   + \gamma_{b}w^{b}_{4}|h^{-}_{4}|(\lambda_{4}^{-}-\dot{y}_{\alpha})
   - w^{m}_{1}|h^{+}_{1}||\lambda_{1}^{+}-\dot{y}_{\alpha}|\big) \\
&+{}(\Ka_{1}A_{W^{+}_{1}} + \widehat{\Ka})|q^{-}_{1}||\lambda_{1}^{-}-\dot{y}_{\alpha}|
- (\Ka_{1}A_{W^{+}_{1}} + \widehat{\Ka})|q^{+}_{1}||\lambda_{1}^{+}-\dot{y}_{\alpha}| \\
& + 2B\Ka_{1}w^{b}_{1}|h^{-}_{1}||\lambda_{1}^{-}-\dot{y}_{\alpha}|
  - 2B\Ka_{1}|q^{-}_{4}|(\lambda_{4}^{-}-\dot{y}_{\alpha}) \\
& + \big(2B\Ka_{1} + \Ka_{1}A_{W^{+}_{4}}
+ \widehat{\Ka}\big)w^{m}_{4}|K_{14}h^{-}_{4}|(\lambda_{4}^{+}-\dot{y}_{\alpha})\\
& -{}\big(2B\Ka_{1} + \Ka_{1}A_{W^{+}_{4}} + \widehat{\Ka}\big)|q^{-}_{4}|(\lambda_{4}^{-}-\dot{y}_{\alpha}) \\
&+\mathcal{O}(1)B\Big(\vartheta + \sum_{i=1,4}|q_{i}^{-}|\Big)\\
=&\, -2(1-\gamma_{b})B\Ka_{1}w^{b}_{4}|h^{-}_{4}|(\lambda_{4}^{-}-\dot{y}_{\alpha})\\
 & + (\Ka_{1}A_{W^{+}_{4}} + \widehat{\Ka})\big(w^{m}_{4}|K_{14}h^{-}_{4}|(\lambda_{4}^{+}-\dot{y}_{\alpha})
   -w^{b}_{4}|h^{-}_{4}|(\lambda_{4}^{-}-\dot{y}_{\alpha})\big)  \\
 &+2B\Ka_{1}w^{b}_{1}(|h^{+}_{1}| + |h^{-}_{1}|)|\lambda_{1}^{-}-\dot{y}_{\alpha}|
   -2B\Ka_{1}w^{m}_{1}|h^{+}_{1}||\lambda_{1}^{+}-\dot{y}_{\alpha}| \\
 &+2B\Ka_{1}w^{m}_{4}|K_{14}h^{-}_{4}|(\lambda_{4}^{+}-\dot{y}_{\alpha})
  -2B\Ka_{1}w^{b}_{4}|h^{-}_{4}|(\lambda_{4}^{-}-\dot{y}_{\alpha}) \\
  &+ (\Ka_{1}A_{W^{+}_{1}} + \widehat{\Ka})w^{b}_1|h^{-}_{1}||\lambda_{1}^{-}-\dot{y}_{\alpha}|
   - (\Ka_{1}A_{W^{+}_{1}} + \widehat{\Ka})w^{m}_1|h^{+}_{1}||\lambda_{1}^{+}-\dot{y}_{\alpha}| \\
  &+{}\mathcal{O}(1)B\Big(\vartheta + \sum_{i=1,4}|q_{i}^{-}|\Big)\\
\leq&\, \mathcal{O}(1)B\vartheta.
\end{align*}

{\it Case $2$}. {\it The weak wave $\alpha$ between the two strong vortex sheets/entropy waves
in $U$ and $V$ is crossed}.

\smallskip
For  $j = 1$, we have
{\small
\begin{align*}
E_{1} &= |q^{\pm}_{1}|(W^{+}_{1} - W^{-}_{1})(\lambda^{\pm}_{1} - \dot{y}_{\alpha}) + W^{\mp}_{1}\left(|q^{+}_{1}|(\lambda^{+}_{1}-\dot{y}_{\alpha})
- |q^{-}_{1}|(\lambda^{-}_{1}-\dot{y}_{\alpha})\right)  \nonumber \\[1mm]
&\le
\begin{cases}
 - \Ka_{1}|q^{\pm}_{1}||\alpha||\lambda^{\pm}_{1} - \dot{y}_{\alpha}|
  +\mathcal{O}(1)\big(B\Ka_{1}+1\big)\big(|q_1^+-q_1^-|+|q_1^-||\alpha|\big)+\mathcal{O}(1)|\alpha|&\\
     &\mbox{when $k_\alpha=2,3,4$,}$\quad$\\
  \mathcal{O}(1)\big(B\Ka_{1}+1\big)\big(|q_1^+-q_1^-|+|q_1^-||\alpha|\big)+\mathcal{O}(1)|\alpha|
    &\mbox{when $k_\alpha=1$}.
\end{cases}
\end{align*}
}

\smallskip
For $j = 2, 3$, we have
\begin{align}
E_{j}
&= B\left(\big(W^{+}_{j} - W^{-}_{j}\big)\big(\lambda^{\pm}_{j} - \dot{y}_{\alpha}\big)
    + W^{\mp}_{j}\big(\lambda^{\pm}_{j}-\lambda^{\mp}_{j}\big)\right) \nonumber \\
&\le
\begin{cases}
-B\left(\Ka_{1}|\alpha| |\lambda^{+}_{j} - \dot{y}_{\alpha}| - \mathcal{O}(1) |\alpha| \right)
\quad &\mbox{when $k_\alpha=1,4$},\\
B\Ka_{1}|\alpha|\big(\vartheta+\mathcal{O}(1)|h_4^-|\big)
   &\mbox{when $k_\alpha=2,3$}.
\end{cases}
\end{align}

For $j = 4$, we have
{\small
\begin{align}
E_{4} &= |q^{\pm}_{4}|(W^{+}_{4} - W^{-}_{4})(\lambda^{\pm}_{4} - \dot{y}_{\alpha})
  + W^{\mp}_{4}\big(|q^{+}_{4}|(\lambda^{+}_{4}-\dot{y}_{\alpha})
  - |q^{-}_{4}|(\lambda^{-}_{4}-\dot{y}_{\alpha})\big)  \nonumber \\[1mm]
&\le
\begin{cases}
-\Ka_{1}|q^{\pm}_{4}||\alpha||\lambda^{\pm}_{4} - \dot{y}_{\alpha}|
+ (2B\Ka_{1}+\mathcal{O}(1))\big((|q^{+}_{4}| - |q^{-}_{4}|)
 (\lambda^{+}_{4} - \dot{y}_{\alpha})&+|q^{-}_{4}|(\lambda_4^+-\lambda_4^-)\big)\\
   &\mbox{when $k_\alpha=1,2,3$},\,\\
 \mathcal{O}(1)(B\Ka_{1}+1)\big(|q^{+}_{4}-q^{-}_{4}| + |q_4^-||\alpha|\big)
  +\mathcal{O}(1)|\alpha|
   &\mbox{when $k_\alpha=4$}.
\end{cases}
\end{align}
}

Then we obtain that, when $k_\alpha=1,4$,
\begin{align*}
\sum_{j=1}^{4}E_{j}
\le&\,  -B\Ka_{1}\big(|\lambda_2^\pm- \dot{y}_{\alpha}| + |\lambda_3^\pm- \dot{y}_{\alpha}|\big)|\alpha|
     +\mathcal{O}(1)(1+B)|\alpha|\\
   &+2B\Ka_1\mathcal{O}(1)
     \big(|q_1^+-q_1^-|+|q_4^+-q_4^-|+(|q_1^-|+ |q_4^-|)|\alpha|\big);
\end{align*}
and, when $k_\alpha=2,3$,
\begin{align*}
\sum_{j=1}^{4}E_{j}
\le&\,  -\Ka_{1}|\alpha|\big(|q_1^\pm||\lambda_1^\pm- \dot{y}_{\alpha}| + |q_4^\pm||\lambda_4^\pm- \dot{y}_{\alpha}|\big)
   +B\Ka_1|\alpha|\big(\vartheta+\mathcal{O}(1)|h_4^-|\big)\\
  &+ \mathcal{O}(1)(B\Ka_1+1)(w_1^m+w_4^m)|h_4^-||\alpha|.
\end{align*}
Choosing $w_j^m, j=1,4$, and $B$ small enough, and $\Ka_1$ large enough, we conclude
$$
\sum_{j=1}^{4}E_{j}\le \mathcal{O}(1)\vartheta.
$$

\smallskip
{\it Case $3$}. {\it The second strong vortex sheet/entropy wave $\alpha$ in $U$ or $V$ is crossed}.
For this case, by Lemma 2.6, we have
\begin{align}
h^{-}_{1} &= K_{21}h^{+}_{1}, \label{5.5} \\
h^{-}_{4} &= h^{+}_{4} + K_{24}h^{+}_{1}. \label{5.6}
\end{align}

\noindent Moreover, the essential estimate $|K_{24}| < 1$ in Lemma 2.6 ensures
the existence of desired weights $w^{a}_{1}$ and $w^{a}_{4}$ in the following manner.
\vspace{3mm}

\noindent
$\textbf{Lemma 5.2.}$
\emph{There exist $w^{a}_{1}$, $w^{a}_{4}$, and $\gamma_{a}$ such that}
\begin{equation}\label{5.7}
\frac{w^{a}_{4}}{w^{a}_{1}}\left|\frac{\lambda_{4}^{+}
- \lambda_{2,3}}{\lambda_{1}^{+} - \lambda_{2,3}}\right|K_{24} < \gamma_{a} < 1.
\end{equation}
\qquad

With Lemma 5.2, we estimate $E_{j}$ for $j= 1, \ldots, 4$ as follows:
By \eqref{5.5},
\begin{align}
E_{1} &= |q_{1}^{-}|(\lambda_{1}^{-}-\dot{y}_{\alpha})(W_{1}^{+} - W_{1}^{-}) + W^{+}_{1}\left(|q^{+}_{1}|(\lambda_{1}^{+}-\dot{y}_{\alpha})
-|q^{-}_{1}|(\lambda^{-}_{1}-\dot{y}_{\alpha})\right)          \nonumber \\
&= -2B\Ka_{1}|q^{-}_{1}||\lambda_{1}^{-}-\dot{y}_{\alpha}|
+ (4B\Ka_{1} + \Ka_{1}A_{W^{+}_{1}} + \widehat{\Ka})|q^{+}_{1}|(\lambda_{1}^{+}-\dot{y}_{\alpha}) \nonumber \\
&\quad +w^{m}_{1}|K_{21}h^{+}_{1}|(4B\Ka_{1} + A_{W^{+}_{1}} + \widehat{\Ka})|\lambda_{1}^{-}-\dot{y}_{\alpha}| \nonumber \\
&= -2B\Ka_{1}|q^{-}_{1}||\lambda_{1}^{-}-\dot{y}_{\alpha}|
+ w^{m}_{1}|K_{21}h^{+}_{1}|(4B\Ka_{1} + \Ka_{1}A_{W^{+}_{1}} + \widehat{\Ka})|\lambda_{1}^{-}-\dot{y}_{\alpha}| \nonumber \\
&\quad -w_{1}^{a}|h^{+}_{1}|(2B\Ka_{1} + \Ka_{1}A_{W^{+}_{1}} + \widehat{\Ka})|\lambda_{1}^{+}-\dot{y}_{\alpha}|
   - 2B\Ka_{1} w_{1}^{a}|h^{+}_{1}||\lambda_{1}^{+}-\dot{y}_{\alpha}|, \nonumber
\end{align}
where $W^{+}_{1} = W_{1}(y_{\alpha}^+) = 4B\Ka_{1} + \Ka_{1}A_{W^{+}_{1}} + \widehat{\Ka}$ with $\widehat{\Ka}$ determined by \eqref{5.5a},
$A_{W^{+}_{1}} = F_{1}(y_{\alpha}^+) + H_{1}(y_{\alpha}^+)$ here is the total strength
of all the weak waves in $U$ and $V$ which approach the $1$-wave $q_{1}^{+} = q_{1}(y_{\alpha}^+)$,
and $4B\Ka_{1}$ is from weight $G_{1}(y_{\alpha}^{+})$.

\smallskip
For $j= 2, 3$, $W^{+}_{j} = W^{-}_{j}$ so that \eqref{4.6} reduces to
\begin{equation*}
E_{j}= W^{+}_{j}\big(|q_{j}^{+}|(\lambda_{j}^{+}-\dot{y}_{\alpha})
- |q_{j}^{-}|(\lambda_{j}^{-}-\dot{y}_{\alpha})\big)
\le \mathcal{O}(1)B\Big(\vartheta
+ \sum_{i=1,4}|q_{i}^{+}|
\Big).
\end{equation*}

\smallskip
By \eqref{5.6}--\eqref{5.7}, we have
\begin{align}
E_{4} &= |q_{4}^{-}|(\lambda_{4}^{-}-\dot{y}_{\alpha})(W_{4}^{+} - W_{4}^{-})
   + W^{+}_{4}\left(|q^{+}_{4}|(\lambda_{4}^{+}-\dot{y}_{\alpha})
-|q^{-}_{4}|(\lambda^{-}_{4}-\dot{y}_{\alpha})\right)          \nonumber \\
&= W^{+}_{4}\left(|q^{+}_{4}|(\lambda_{4}^{+}-\dot{y}_{\alpha})
-|q^{-}_{4}|(\lambda^{-}_{4}-\dot{y}_{\alpha})\right)          \nonumber \\
&\le W^{+}_{4}\left(w^{a}_{4}|h^{-}_{4}|(\lambda_{4}^{+}-\dot{y}_{\alpha})
+ w^{a}_{4}K_{24}|h^{+}_{1}|(\lambda_{4}^{+}-\dot{y}_{\alpha})
- w^{m}_{4}|h^{-}_{4}|(\lambda_{4}^{-}-\dot{y}_{\alpha})\right) \nonumber \\
&\le 2B\Ka_{1}\left(w^{a}_{4}|h^{-}_{4}|(\lambda_{4}^{+}-\dot{y}_{\alpha})
  + \gamma_{a}w^{a}_{1}|h^{+}_{1}||\lambda_{1}^{+}-\dot{y}_{\alpha}|
- w^{m}_{4}|h^{-}_{4}|(\lambda_{4}^{-}-\dot{y}_{\alpha})\right) \nonumber \\
&\quad + (\Ka_{1}A_{W^{+}_{4}} + \widehat{\Ka})|q^{+}_{4}|(\lambda_{4}^{+}-\dot{y}_{\alpha})
  - (\Ka_{1}A_{W^{+}_{4}} + \widehat{\Ka})|q^{-}_{4}|(\lambda_{4}^{-}-\dot{y}_{\alpha}), \nonumber
\end{align}
where $W^{+}_{4} = W_{4}(y_{\alpha}^+) = 2\Ka_{1}B + \Ka_{1}A_{W^{+}_{4}} + \widehat{\Ka}$ with $\widehat{\Ka}>0$ determined
by \eqref{5.5a},
$A_{W^{+}_{4}} = F_{4}(y_{\alpha}^+) + H_{4}(y_{\alpha}^+)$ here is the total strength
of all the weak waves in $U$ and $V$ which approach the $4$-wave $q_{4}^{+} = q_{4}(y_{\alpha}^+)$,
and $2B\Ka_{1}$ is from weight $G_{4}(y_{\alpha}^{+})$.

\smallskip
For the weighted $L^{1}$--strength $q_{j}(y)$ in (4.1),
when $w^{a}_{4}$ is small enough relative to $w^{m}_{4}$, $w^{a}_{1}$ is large enough
relative to $w^{m}_{1}$, $\Ka_{1}$ is large enough, applying \eqref{5.5}--\eqref{5.6},
suitably small total variation of $U$ and $V$ yields
\begin{align*}
\sum^{4}_{j=1}E_{j}
\leq&\, 2B\Ka_{1}\left(w^{a}_{4}|h^{-}_{4}|(\lambda_{4}^{+}-\dot{y}_{\alpha}) + \gamma_{a}w^{a}_{1}|h^{+}_{1}||\lambda_{1}^{+}-\dot{y}_{\alpha}|
  - w^{m}_{4}|h^{-}_{4}|(\lambda_{4}^{-}-\dot{y}_{\alpha})\right) \\
 &+{}(\Ka_{1}A_{W^{+}_{4}} + \widehat{\Ka})|q^{+}_{4}|(\lambda_{4}^{+}-\dot{y}_{\alpha})- (\Ka_{1}A_{W^{+}_{4}} + \widehat{\Ka})|q^{-}_{4}|(\lambda_{4}^{-}-\dot{y}_{\alpha}) \\
 &-2B\Ka_{1}|q^{-}_{1}||\lambda_{1}^{-}-\dot{y}_{\alpha}|
  - 2B\Ka_{1}|q^{+}_{1}||\lambda_{1}^{+}-\dot{y}_{\alpha}| \\
 &+ \big(4B\Ka_{1} + \Ka_{1}A_{W^{+}_{1}} + \widehat{\Ka}\big)w^{m}_{1}|K_{21}h^{+}_{1}||\lambda_{1}^{-}-\dot{y}_{\alpha}|\\
 &-{}\big(2B\Ka_{1} + \Ka_{1}A_{W^{+}_{1}} + \widehat{\Ka}\big)w^{a}_{1}|h^{+}_{1}||\lambda_{1}^{+}-\dot{y}_{\alpha}| \\
  &+\mathcal{O}(1)B\Big(\vartheta + \sum_{i=1,4}|q_{i}^{+}|
  \Big)\\
=&-2(1-\gamma_{a})B\Ka_{1}w^{a}_{1}|h^{+}_{1}||\lambda_{1}^{+}-\dot{y}_{\alpha}| \\
& + (\Ka_{1}A_{W^{+}_{1}} + \widehat{\Ka})\big(w^{m}_{1}|K_{21}h^{+}_{1}||\lambda_{1}^{-}-\dot{y}_{\alpha}|-w^{a}_{1}|h^{+}_{1}||\lambda_{1}^{+}-\dot{y}_{\alpha}|\big)  \\
&+2B\Ka_{1}w^{a}_{4}|h^{-}_{4}|(\lambda_{4}^{+}-\dot{y}_{\alpha})
-2B\Ka_{1}w^{m}_{4}|h^{-}_{4}|(\lambda_{4}^{-}-\dot{y}_{\alpha}) \\
&+2B\Ka_{1}(-w^{m}_{1}|K_{21}h^{+}_{1}| + 2w^{m}_{1}|K_{21}h^{+}_{1}|)|\lambda_{1}^{-}-\dot{y}_{\alpha}| - 2B\Ka_{1}w^{a}_1|h^{+}_{1}||\lambda_{1}^{+}-\dot{y}_{\alpha}| \\
&+ (\Ka_{1}A_{W^{+}_{4}} + \widehat{\Ka})w^{a}_{4}|h^{+}_{4}|(\lambda_{4}^{+}-\dot{y}_{\alpha}) - (\Ka_{1}A_{W^{+}_{4}} + \widehat{\Ka})w^{m}_{4}|h^{-}_{4}|(\lambda_{4}^{-}-\dot{y}_{\alpha}) \\
&+{}\mathcal{O}(1)B\Big(\vartheta + \sum_{i=1,4}|q_{i}^{+}|
\Big) \\
\leq&\, 0,
\end{align*}
which implies \eqref{4.10}.

\smallskip
{\it Case $4$}. {\it Close to the Lipschitz wall boundary}.
This case differs from the Cauchy problem. Here we will use the particular property
of the boundary condition \eqref{1.9}: The flows of $U$ and $V$ are tangent to the Lipschitz wall,
which implies that they must be parallel with each other along the boundary.
Then a piecewise constant weak solution is constructed only along the Hugoniot curves
determined by the Riemann data $U(b)$ and $V(b)$,
the states of solutions $U$ and $V$, respectively, close to the boundary.
\vspace{4mm}

\noindent
\textbf{Lemma 5.3.}
\emph{Let $U(b) =(\breve{u}, \breve{v}, \breve{p}, \breve{\rho})$
and $V(b) = (\tilde{u}, \tilde{v}, \tilde{p}, \tilde{\rho})$
be two states in a small neighborhood $O_{\varepsilon}(U_{-})$ of $U_{-}$ satisfying
$\frac{\breve{v}}{\breve{u}} =\frac{\tilde{v}}{\tilde{u}} = \dot{z}_{b}$
and $\breve{v}, \tilde{v} \approx 0$.
Denote by $h_{j}(b)$ as the strength of the $j^{th}$-shock in the Riemann problem
determined by $U(b)$ and $V(b)$,
and denote by $\lambda_{j}$ as the corresponding $j^{th}$-characteristic speed. Then}
\begin{eqnarray}
&&|\lambda_{j} - \dot{z}_{b}| \sim |h_{1}(b)| \qquad \mbox{for $j=2, 3$}, \\[2mm]
&&|h_{4}(b)| \leq |h_{1}(b)|
 + \mathcal{O}(1)|h_{2}(b)||\lambda_{2} - \dot{z}_{b}| + |h_{1}(b)|\mathcal{O}(1)|\dot{z}_{b}|, \\[2mm]
&&|h_{1}(b)| = \bar{\mathcal{O}}(1)|h_{4}(b)| \qquad \mbox{with $\frac{1}{2} < \bar{\mathcal{O}}(1) < \frac{3}{2}$},
\end{eqnarray}
\noindent \emph{where $\dot{z}_{b}$ is the slope of the Lipschitz wall.}
\vspace{3mm}

\textbf{Proof.} We prove this by analyzing the two cases.

\smallskip
\emph{Case $1$. $h_{1}(b) =0 $ and $h_{4}(b)=0$ that corresponds to the case $\breve{p} = \tilde{p}$}.
Starting at state $U_{b}$, we move along the Hugoniot curves of the second and third families
to reach $V_{b}$.
Note that these two families are the contact Hugoniot curves,
so that $\lambda_{2}$ and $\lambda_{3}$ are constant along the Hugoniot curves.
Given that $\lambda_{2,3} =\frac{v}{u}$,
$\textbf{r}_{2} = (1, \frac{v}{u}, 0, 0)^{\top}$,
and $\textbf{r}_{3} = (0, 0, 0, \rho)^{\top}$,
$\frac{v}{u}$ remains unchanged as the initial value $\frac{v(U_{b})}{u(U_{b})}$,
{\it i.e.}, $\dot{z}_{b}$ in this process by the boundary condition \eqref{1.9}.
Therefore, we conclude that $\lambda_{2, 3} = \dot{z}_{b}$, equivalently,
$$
\dot{z}_{b} - \lambda_{2, 3} = 0.
$$

\smallskip
\emph{Case $2$. $h_{1}(b) \neq 0$ that corresponds to $\breve{p} \neq \tilde{p}$}.
Starting at state $U(b)$, we move along the 1-Hugoniot curve to reach $U_{1}$,
then possibly move along the 2-contact Hugoniot curve to reach $U_{2}$,
the 3-Hugoniot curve to reach $U_{3}$, and the 4-Hugoniot curve to reach $V(b)$.

To make clear some essential relations among the strengths: $h_{1}(b), h_{2}(b), h_{3}(b), \text{and } h_{4}(b)$,
we project $(u, v, p, \rho)$ onto the $(u, v)$--plane.
Denote $\textbf{r}_{1}|_{u}$ as the projection of $\textbf{r}_{1}$ onto the $u$-axis,
$\textbf{r}_{2}|_{(u,v)}$ as the projection of $\textbf{r}_{2}$ onto the $(u,v)$--plane, {\it etc}.
Then, at the background state $U_{-}$,
$$
\textbf{r}_{1}|_{u} = -\textbf{r}_{4}|_{u}, \hspace{2mm}
\textbf{r}_{1}|_{v} = \textbf{r}_{4}|_{v}, \hspace{2mm}
\textbf{r}_{1}|_{(p, \rho)} = -\textbf{r}_{4}|_{(p, \rho)}, \hspace{2mm}
\textbf{r}_{2} = \textbf{r}_{2}|_{(u,v)}, \hspace{2mm} \textbf{r}_{3}|_{(u,v)} = 0.
$$

\noindent We first note that $h_{4}(b) \neq 0$.
Given that $\textbf{r}_{1}|_{(u,v)} = k_{1}(-\lambda_{1}, 1)^{\top}$ along with finite characteristic
speeds $\lambda_{1}$ and $\dot{z}_{b} \approx 0$, then
$ \dot{z}_{b} < -\frac{1}{\lambda_{1}}$ near state $U_{-}$.
Thus, we can conclude that, in the $(u,v)$--plane,
$\frac{{\rm d}v}{{\rm d}u}$ along the 1-curve is always larger than $\dot{z}_{b}$.
This implies that $\frac{v(U_{1})}{u(U_{1})} \neq \frac{v(U_{b})}{u(U_{b})}$.
Meanwhile,
$\frac{v(U_{1})}{u(U_{1})} = \frac{v(U_{2})}{u(U_{2})}
=\frac{v(U_{3})}{u(U_{3})}$ and $\frac{v(V_{b})}{u(V_{b})} =\frac{v(U_{b})}{u(U_{b})}$.
This implies that
$$
\frac{v(U_{1})}{u(U_{1})} = \frac{v(U_{2})}{u(U_{2})} = \frac{v(U_{3})}{u(U_{3})} \neq \frac{v(V_{b})}{u(V_{b})}.
$$
Therefore, we conclude that there is some distance
along the 4-Hugoniot curve to reach $V_{b}$ so that
$h_{4} \neq 0$.

Next, we present an essential estimate to  bound $\left|h_{1}\right|$ more precisely
in terms of $\left|h_{4}\right|$.
To that end, define the signed length of $(U_{1} - U_{b})|_{(u,v)}$
and $(V_{b} - U_{3})|_{(u,v)}$ by $d_{1}$ and $d_{4}$ on the $(u,v)$--plane:
\begin{equation}
 d_{1} =
 \begin{cases}
 \|(U_{1} - U_{b})|_{(u,v)}\|&\mbox{ if $h_{1} > 0$,} \\[2mm]
 -\|(U_{1} - U_{b})|_{(u,v)}\|&\mbox{ if $h_{1} < 0$,} \nonumber \\
 \end{cases}
\end{equation}
and
\begin{equation}
 d_{4} = \begin{cases}
 \|(V_{b} - U_{3})|_{(u,v)}\|&\mbox{ if $h_{4} > 0$,} \\[1mm]
 -\|(V_{b} - U_{3})|_{(u,v)}\|&\mbox{ if $h_{4} < 0$.} \nonumber \\
        \end{cases}
\end{equation}

Note that
$$
|\lambda_{2} - \dot{z}_{b}| = \mathcal{O}(1)|d_{1}| = \mathcal{O}(1)|h_{1}(b)|.
$$
Since $\lambda_{2} = \frac{v(U_1)}{u(U_1)} = \frac{v(U_2)}{u(U_2)} = \lambda_{3}$,
we can similarly conclude
$$
|\lambda_{3} - \dot{z}_{b}| = \mathcal{O}(1)|d_{1}| = \mathcal{O}(1)|h_{1}(b)|,
$$
by using the following projections on the $(u,v)$--plane:
$$
\textbf{r}_{1}|_{u} = -\textbf{r}_{4}|_{u},\quad  \textbf{r}_{1}|_{v} = \textbf{r}_{4}|_{v},
\quad \textbf{r}_{2} = \textbf{r}_{2}|_{(u,v)}, \quad \textbf{r}_{3}|_{(u,v)} = 0.
$$

Moreover, we note that
$$
-d_{4} = \mathcal{O}(1) h_{2}(b)(\lambda_{2} - \dot{z}_{b}) + \tilde{d},
$$
where $\tilde{d}\cos \varphi_{1} = d_{1}\cos \varphi_{2}$,
$\varphi_{1}$ denotes the angle between $(1, \dot{z}_{b})$ and $\textbf{r}_{4}|_{(u,v)}$,
$\varphi_{2}$ denotes the angle between $\textbf{r}_{1}|_{(u,v)}$ and $(1, \dot{z}_{b})$,
$\varphi_{1} = \varphi_{2} + 2\alpha$ for $\alpha= {\rm arctan}(\dot{z}_{b})$, and
\begin{align}
\tilde{d} &= d_{1}\frac{\cos \varphi_{2}}{\cos \varphi_{1}}
=d_{1}\frac{\cos (\varphi_{1}-2\alpha)}{\cos \varphi_{1}}
=d_{1}\frac{\cos \varphi_{1}\cos(2\alpha)+\sin \varphi_{1} \sin(2\alpha)}{\cos \varphi_{1}} \nonumber \\
&= d_{1}\big(\text{cos}(2\alpha) + \mathcal{O}(1)\text{sin}(2\alpha)\big)
 = d_{1}\big(1 + \mathcal{O}(1)\alpha\big)
 = d_{1}\big(1 + \mathcal{O}(1)\dot{z}_{b}\big), \nonumber
\end{align}
so that
\begin{equation*}
-d_{4} = \mathcal{O}(1)h_{2}(b)(\lambda_{2} - \dot{z}_{b}) + d_{1}\big(1 + \mathcal{O}(1)\dot{z}_{b}\big).
\end{equation*}

At $U_{-}$,  $\textbf{r}_{1}|_{(u, p, \rho)}$ = {-}$\textbf{r}_{4}|_{(u, p, \rho)}$
and $\textbf{r}_{1}|_{v} = \textbf{r}_{4}|_{v}$, which implies
$$
\frac{d_{1}}{h_{1}} = \frac{d_{4}}{h_{4}}.
$$
Thus, we obtain the following key estimate:
\begin{equation}\label{5.7a}
-h_{4}(b) = \mathcal{O}(1)h_{2}(b)(\lambda_{2} - \dot{z}_{b}) + h_{1}(b)\big(1 + \mathcal{O}(1)\dot{z}_{b}\big).
\end{equation}
Estimate \eqref{5.7a} now implies
\begin{align}
|h_{4}(b)| &\leq |h_{1}(b)| + \mathcal{O}(1)|h_{2}(b)||\lambda_{2} - \dot{z}_{b}| + |h_{1}(b)|\mathcal{O}(1)|\dot{z}_{b}|
\nonumber \\
&\leq |h_{1}(b)| + \mathcal{O}(1)\big(|h_{2}(b)| + |\dot{z}_{b}|\big)|h_{1}(b)|, \nonumber
\end{align}
which yields
$$
|h_{1}(b)| = \bar{\mathcal{O}}(1)|h_{4}(b)|
\hspace{7mm} \text{with $\frac{1}{2} < \bar{\mathcal{O}}(1) < \frac{3}{2}$},
$$
given that $|h_{2}(b)| + |\dot{z}_{b}|$ is always small enough.
This is guaranteed by the sufficiently small total variation of the initial
perturbation $\widetilde{U_{0}}$ and the boundary perturbation. This completes the proof.

\vspace{5mm}
Notice that the requirement $\frac{\bar{v}}{\bar{u}} = \frac{\hat{v}}{\hat{u}} = \dot{z}_{b}$
in Lemma 5.3 is just the boundary condition \eqref{1.9} because $\dot{z}_{b}$ here is the slope of the Lipschitz wall.

\vspace{3mm}
Applying Lemma 5.3 now yields
\begin{align}
E_{b,1} &= |q_{1}(b)|W_{1}(b)(-\dot{z}_{b} + \lambda_{1})
\nonumber \\&
= -4B\Ka_{1} w_{1}^{b}|h_{1}(b)|\, |\lambda_{1}| + \mathcal{O}(1)|h_{1}(b)|
 \nonumber \\&
=-4B \Ka_1 w_{1}^{b}|h_{1}(b)||\lambda_{1}| + \mathcal{O}(1)|h_{4}(b)|,  \nonumber \\[2mm]
E_{b,j} &= |q_{j}(b)|W_{j}(b)(-\dot{z}_{b} + \lambda_{j})
 = \mathcal{O}(1) w^{b}_{j}|h_{j}(b)|(-\dot{z}_{b}+ \lambda_{j})
 = \mathcal{O}(1)h_{1}(b) = \mathcal{O}(1)h_{4}(b) \hspace{4mm} \mbox{for $j = 2, 3$}, \nonumber \\[2mm]
E_{b,4} &= |q_{4}(b)|W_{4}(b)(-\dot{z}_{b} + \lambda_{4})\nonumber\\
 &= 4B\Ka_1 w^{b}_{4}|h_{4}(b)||\lambda_{1}| + \mathcal{O}(1)|h_{4}(b)| \nonumber \\
&\leq 4B\Ka_{1}|\lambda_{1}|w^{b}_{4}\big(|h_{1}(b)| + \mathcal{O}(1)|h_{2}(b)||\lambda_{2} - \dot{z}_{b}| + \mathcal{O}(1)|h_{1}(b)||\dot{z}_{b}|\big) + \mathcal{O}(1)|h_{4}(b)|. \nonumber
\end{align}
Using Lemma 5.1, we can choose $w^{b}_{1}$ and $w^{b}_{4}$ such that
$$
w^{b}_{4} < w^{b}_{1}.
$$
Then, with the total variation of the incoming flow perturbation and the boundary perturbation
small enough and $\Ka_{1}$ large enough, we have
\begin{align*}
\sum_{j=1}^{4}E_{b,j}
&= 4B\Ka_1(w^{b}_{4} - w^{b}_{1})|h_{1}(b)||\lambda_{1}|
+ \mathcal{O}(1)B\Ka_{1} |\lambda_{1}| h^{b}_{4}\big(|h_{2}(b)| + |\dot{z}_{b}|\big)|h_{1}(b)|
 +\mathcal{O}(1)|h_{4}(b)|\\
&\leq \bar{\mathcal{O}}(1)4B\Ka_{1} (w^{b}_{4} - w^{b}_{1})|h_{4}(b)|\, |\lambda_{1}|
+ \mathcal{O}(1)B\Ka_{1}|\lambda_{1}| w^{b}_{4}\big(|h_{2}(b)| + |\dot{z}_{b}|\big)|h_{1}(b)|
+ \mathcal{O}(1)|h_{4}(b)| \leq 0,
\end{align*}
provided that $|h_{2}(b)| + |\dot{z}_{b}|$ is sufficiently small.
This is guaranteed since the total variation of the incoming flow perturbation
and the boundary perturbation are sufficiently small.

\section{Existence of A Semigroup of Solutions}

As a corollary of the essential estimates in \S 3--\S 5, we can now establish both the existence
of semigroup $\mathscr{S}$ generated by the wave-front tracking method
and the Lipschitz continuity of $\mathscr{S}$.
\vspace{4mm}

\noindent \textbf{Lemma 6.1.}
\emph{If $TV(\widetilde{U_0}(\cdot))+ TV(g^{\prime}(\cdot))$ is
small enough, then the map{\rm :}
$$
(\overline{U}(\cdot), x) \mapsto U^{\vartheta}(x, \cdot)\mathrel{\mathop:}= \mathscr{S}^{\vartheta}_{x}(\overline{U}(\cdot))
$$
produced by the wave-front tracking algorithm is a uniformly Lipschitz continuous semigroup satisfying the properties}:
\begin{itemize}
\item[(\rmnum{1})] $\mathscr{S}^{\vartheta}_{0}\overline{U} = \overline{U}$, and
 $\mathscr{S}^{\vartheta}_{x_1}\mathscr{S}^{\vartheta}_{x_2}\overline{U}
    = \mathscr{S}^{\vartheta}_{x_1 + x_2}\overline{U}$ $\quad$ for all $x_1, x_2 \geq 0$;

\smallskip
\item[(\rmnum{2})] $\norm{\mathscr{S}^{\vartheta}_{x}\overline{U} - \mathscr{S}^{\vartheta}_{x}\overline{V}}_{L^1} \leq C\norm{\overline{U} - \overline{V}}_{L^1} + C\vartheta x$ $\quad$ for all $x \ge 0$.
\end{itemize}

\smallskip
\noindent \textbf{Proof.} Since $\mathscr{S}^{\vartheta}$ is generated
by the wave-front tracking algorithm,
property (i) is immediate.
Next, property (ii) is proved as follows:
Take a pair of front tracking $\vartheta$-approximate solutions $U^{\vartheta}$ and $V^{\vartheta}$
of problem \eqref{2.1}--\eqref{1.8}
with $\overline{U}(\cdot)$ and $\overline{V}(\cdot)$ as the initial data, respectively.
Using \eqref{4.5} and \eqref{4.14a}, at any $x \geq 0$, we have
\begin{align}
\norm{U^{\vartheta}(x) - V^{\vartheta}(x)}_{L^1}
\leq {C}\Phi(U^{\vartheta}(x), V^{\vartheta}(x))
\leq {C}\Phi(U^{\vartheta}(0), V^{\vartheta}(0)) + C\vartheta x
\leq C\norm{\overline{U} - \overline{V}}_{L^1} + C\vartheta x.
\end{align}
Therefore, the $\vartheta$-semigroup is Lipschitz continuous.

\smallskip
For a given $\nu_{0} > 0$, we define the domain:
\begin{equation*}
\mathcal{D} =  \mbox{cl}
\left\{
U\,:\, \mathbb{R} \mapsto \mathbb{R}^4\,\, \left|\begin{array}{ll}
&\mbox{$\exists$ one point $y^{i} \in \mathbb{R}$ and $U_\pm$ such that}\\
&U - \tilde{U} \in L^1(\mathbb{R}; \mathbb{R}^4)
\mbox{ and } {\rm TV}(U - \tilde{U}) \leq \nu_{0}\\
&\mbox{with } \tilde{U}(y) =
\begin{cases} U_{-} \qquad &\mbox{for $g(x) \leq y \leq y^{i}$}, \\
U_{+} \qquad & \mbox{for $y^{i} < y$}
\end{cases}
\end{array}\right.
\right\}.
\end{equation*}
Given a solution $U(x,y)$ to the initial-boundary value problem \eqref{2.1}--\eqref{1.8},
we note that,
if $U^{x}(y) \mathrel{\mathop:}= U(x,y) \in \mathcal{D}$ at any fixed $x \geq 0$,
then $y^i > g(0) = 0$ at $x = 0$ and $y^i > g(x)$ for $x > 0$ as a strong vortex sheet/entropy wave is present.
\vspace{3mm}

The semigroup $\mathscr{S}$ generated by the wave-front tracking algorithm is provided by the following theorem.

\smallskip
\noindent \textbf{Theorem 6.1.} \emph{If $TV(\widetilde{U_0}(\cdot))+TV(g^{\prime}(\cdot))$
is small enough, then $\mathscr{S}^{\vartheta}$ produced by the wave-front tracking
algorithm is a Cauchy sequence in the $L^1$--norm, so that $\mathscr{S}^{\vartheta}$ converges to a unique limit $\mathscr{S}$ satisfying that
$\mathscr{S}_{x}(\overline{U}) = \lim_{\vartheta \to 0}\mathscr{S}^{\vartheta}_{x}(\overline{U})$ for any $x > 0$.
Then the map $\mathscr{S}: [0, \infty) \times \mathcal{D} \mapsto \mathcal{D}$ is a uniformly Lipschitz semigroup in $L^1$.}
\emph{In particular, the entropy solution to the initial-boundary value problem \eqref{2.1}--\eqref{1.8}
constructed by the wave-front tracking algorithm is unique and $L^1$--stable}.
\vspace{3mm}

Based on the essential estimates in \S 3--\S 5,
Theorem 6.1 can be proved in the same way as the argument given in \ci{Bressan-Colombo-1995}.
Also see Chen-Li \ci{Chen-Li-2008}.

\section{Uniqueness of Entropy Solutions in the Class of Viscosity Solutions}

In this section,
as an immediate consequence of the estimates obtained in \S 4--\S 6,
we find that the semigroup $\mathscr{S}$ produced by the wave-front tracking method is the
only {\it standard Riemann semigroup} (SRS) in the sense of Definition 7.1 given below.
In other words, the semigroup defined by the wave-front tracking method
is the canonical trajectory of the standard Riemann semigroup (SRS).
This yields the uniqueness of entropy solutions in a broader class
of viscosity solutions as introduced by Bressan in \ci{Bressan-1995}.
Furthermore, it coincides with the semigroup
trajectory generated by the wave-front tracking method.

\medskip
the initial-boundary value problem \eqref{2.1}--\eqref{1.8},
is said to have a standard Riemann semigroup
if,
for some small $\nu_{0}$, there exist both a continuous mapping $\mathscr{R}: [0, \infty) \times \mathcal{D} \mapsto \mathcal{D}$
and a constant $L$ satisfying the following properties:
\begin{itemize}
\item[(\rmnum{1})] {\it Semigroup property}: $\mathscr{R}_{0}\overline{U} = \overline{U}$ and
 $\mathscr{R}_{x_1}\mathscr{R}_{x_2}\overline{U} = \mathscr{R}_{x_1 + x_2}\overline{U}$;

\vspace{2mm}
\item[(\rmnum{2})] {\it Lipschitz continuity}:
$\norm{\mathscr{R}_{x}\overline{U} - \mathscr{R}_{x}\overline{V}}_{L^1} \leq L\norm{\overline{U} - \overline{V}}_{L^1}$;

\vspace{2mm}
\item[(\rmnum{3})] {\it Consistency with the Riemann solver}: \emph{Given piecewise constant initial data
$\overline{U} \in \mathcal{D}$, then, for all $x \in [0,\nu_0]$, $U(x, \cdot) = \mathscr{R}_{x}\overline{U}$ coincides
with the solution of problem \eqref{2.1}--\eqref{1.8}
obtained by piecing together the standard Riemann solutions
and the lateral Riemann solutions.}
\end{itemize}

\medskip
Following the argument in \ci{Bressan-1995}, we employ the estimates obtained in \S 4--\S 6 to conclude

\medskip
\noindent \textbf{Theorem 7.1.} \emph{Suppose that problem \eqref{2.1}--\eqref{1.8}
has a standard Riemann
semigroup $\mathscr{R} : [0, \infty) \times \mathcal{D} \mapsto \mathcal{D}$.
Consider the semigroup $\mathscr{S}$ produced by the wave-front tracking method{\rm :}
$\mathscr{S}_{x}(\overline{U}) = \lim_{\vartheta \rightarrow 0}\mathscr{S}^{\vartheta}_{x}(\overline{U})$.
Assume $\overline{U} \in \mathcal{D}$.
Then, for all $x > 0$, $\mathscr{R}_{x}\overline{U} = \mathscr{S}_{x}\overline{U}$}.
\emph{Furthermore, a continuous map $U :[0, X] \mapsto \mathcal{D}$ is a viscosity solution of
problem \eqref{2.1}--\eqref{1.8}
defined in \ci{Bressan-1995}
if and only if
\begin{equation}\label{7.3}
U(x, \cdot) = \mathscr{R}_{x}\overline{U} \hspace{7mm} \textit{for any $x\in [0, T]$}.
\end{equation}
In particular, a continuous map $U :[0, X] \mapsto \mathcal{D}$ is a viscosity solution if and only if
\begin{equation}
U(x, \cdot) = \mathscr{S}_{x}\overline{U} \hspace{7mm} \textit{for any $x\in [0, T]$}.
\end{equation}
}

The proof here follows a similar argument to the one presented in \ci{Bressan-1995}.
The only difference is the strong vortex sheet/entropy wave in our problem.
Nonetheless, one can proceed with the proof by considering the convergence
of the wave-front tracking method which is shown in \S 3.
\vspace{4mm}

\noindent \textbf{Remark 7.1.} In the simpler cases of the isentropic or isothermal Euler flow,
as well as the potential flow, as far as the $L^1$--stability problem is concerned,
we realize the same results as those for the full Euler system \eqref{1.1}.

\bigskip

\noindent
{\bf Acknowledgements:}
The authors would like to thank Yun Pu for his helpful suggestions.
The research of Gui-Qiang G. Chen was supported in part by
the UK Engineering and Physical Sciences Research Council Awards
EP/L015811/1, EP/V008854/1, and EP/V051121/1.
The research of Vaibhav Kukreja was supported in part
by the National Science
Foundation under Grants
DMS-0935967 and DMS-0807551, the UK EPSRC Science and Innovation
Award to the Oxford Centre for Nonlinear PDE (EP/E035027/1).


\bigskip

\begin{thebibliography}{99}

\bib{Baiti-Jenssen-1998}
P. Baiti and K. Jenssen, On the front-tracking algorithm, J. Math. Anal. Appl. 217 (1998), 395--404.


\bib{Bressan-1992}
A. Bressan, Global solutions of systems of conservation laws by wave-front tracking, J. Math. Anal. Appl. 170 (1992), 414--432.

\bib{Bressan-1995}
A. Bressan, The unique limit of the Glimm scheme, Arch. Ration. Mech. Anal. 130 (1995), 205--230.

\bib{Bressan-2000}
A. Bressan, Hyperbolic Systems of Conservations Laws: The One-Dimensional Cauchy Problem, Oxford University Press: Oxford, 2000.

\bib{Bressan-Colombo-1995}
A. Bressan and R.~M. Colombo, The semigroup of $2 \times 2$ conservation laws, Indiana Univ. Math. J. 44 (1995), 677--725.


\bib{Bressan-Liu-Yang-1999}
A. Bressan, T.-P. Liu, and T. Yang, $L^{1}$ stability estimates for $n \times n$ conservation laws, Arch. Ration. Mech. Anal. 149 (1999), 1--22.

\bib{Chen-Feldman2018}
G.-Q. Chen and M. Feldman,
Mathematics of Shock Reflection-Diffraction and von Neumann's Conjecture,
Research Monograph, Annals of Mathematics Studies, 197,
Princeton University Press: Princeton, 2018.

\bib{Chen-Li-2008}
G.-Q. Chen and T.-H. Li,
Well-posedness for two-dimensional steady supersonic Euler flows past a Lipschitz wedge,
J. Diff. Equ. 244 (2008), 1521--1550.

\bibitem{CSV}
G.-Q. Chen, H. Shahgholian, and J.-V. V\'{a}zquez,
Free boundary problems: The forefront of current and future developments,
In: Free Boundary Problems and Related Topics.
Theme Volume: Phil. Trans. R. Soc. A. 373 (2015), 20140285.


\bibitem{CW12}
G.-Q. Chen and Y.-G. Wang,
Characteristic discontinuities and free boundary problems for hyperbolic conservation laws.
In: Nonlinear Partial Differential Equations, 53--81, Abel Symp. 7,
Springer, Heidelberg, 2012.

\bib{Chen-Zhang-Zhu-2007}
G.-Q. Chen, Y.-Q. Zhang, and D.-W. Zhu,
Stability of compressible vortex sheets in steady supersonic Euler flows over Lipschitz walls,
SIAM J. Math. Anal. 38 (2007), 1660--1693.

\bibitem{CDK}
E. Chiodaroli, C. De Lellis, and O. Kreml,
Global ill-posedness of the isentropic system of gas dynamics,
\textit{Comm. Pure Appl. Math.} \textbf{68} (2015), 1157--1190.


\bib{CS1}
J. F. Coulombel and P. Secchi, The stability of compressible vortex sheets in two space
dimensions, Indiana Univ. Math. J. 53 (2004), 941--1012.

\bib{CS2}
J. F. Coulombel and P. Secchi, Nonlinear compressible vortex sheets in two space dimensions,
Ann. Sci. Ec. Norm. Super. 41 (2008), 85--139.


\bib{Courant-Friedrichs-1948}
R. Courant and K.~O. Friedrichs,
Supersonic Flow and Shock Waves, Interscience: New York, 1948.

\bib{Corli-Sable-1997}
A. Corli and M. Sabl$\acute{\text{e}}$-Tougeron,
Stability of contact discontinuities under perturbations of bounded variation, Rend. Sem. Mat. Univ. Podova, 97 (1997), 35--60.

\bib{Dafermos-1972}
C.~M. Dafermos, Polygonal approximations of solutions of the initial value problem for a conservation law, J. Math. Anal. Appl. 38 (1972), 33--41.

\bib{Dafermos-2005}
C.~M. Dafermos, Hyperbolic Conservation Laws in Continuum Physics, Second Edition, Springer-Verlag: Berlin, 2005.

\bib{DiPerna-1976}
R.~J. DiPerna, Global existence of solutions to nonlinear hyperbolic systems of conservation laws, J. Diff. Equ. 20 (1976), 187--212.

\bib{Glimm-1965}
J. Glimm, Solution in the large for nonlinear systems of conservation laws, Comm. Pure Appl. Math. 18 (1965), 697--715.

\bib{Holden-Riesbro-2002}
H. Holden and N. Risebro, Front Tracking for Hyperbolic Conservation Laws, Springer-Verlag: New York, 2002.

\bibitem{KKMM}
C. Klingenberg, O. Kreml, V. M\'{a}cha, and S. Markfelder,
Shocks make the Riemann problem for the full Euler system
in multiple space dimensions ill-posed,
\newblock {\itshape Nonlinearity}, \textbf{33} (2020), 6517--6540.




\bib{Lax}
P.~D. Lax, Hyperbolic Systems of Conservation Laws and the Mathematical Theory of Shock Waves.
CBMS-RCSAM, No. 11. SIAM: Philadelphia, Pa., 1973.

\bib{LEF}
Ph. LeFloch, Hyperbolic Systems of Conservation Laws: The Theory of Classical and Nonclassical Shock Waves, Birkh$\ddot{\text{a}}$user-Verlag: Basel, 2002.

\bib{Lewicka-2000}
M. Lewicka, $L^{1}$ stability of patterns of non-interacting large shock waves,
Indiana Univ. Math. J. 49 (2000), 1515--1537.

\bib{Lewicka-2001}
M. Lewicka, Stability conditions for patterns of noninteracting large shock waves,
SIAM J. Math. Anal. 32 (2001), 1094--1116.

\bib{Lewicka-Trivisa-2002}
M. Lewicka and K. Trivisa, On the $L^{1}$ well posedness of systems of conservation laws near solutions containing two large shocks,
J. Diff. Equ. 179 (2002), 133--177.

\bib{Liu-1997}
T.-P. Liu, The deterministic version of the Glimm scheme, Commun. Math. Phys. 57 (1977), 135--148.

\bib{Liu-Yang-1999}
T.-P. Liu and T. Yang, Well-posedness theory for hyperbolic conservation laws, Comm. Pure Appl. Math. 52 (1999), 1553--1586.

\bib{Sable-1993}
M. Sabl$\acute{\text{e}}$-Tougeron, M$\acute{\text{e}}$thode de Glimm et probl$\grave{\text{e}}$me mixte, Ann. Inst. H. Poincar$\acute{\text{e}}$ Anal. Nonlin$\acute{\text{e}}$aire, 10 (1993), 423--443.
\end{thebibliography}
\end{document}